\documentclass[a4paper,10pt]{amsart}

\usepackage{amsmath,amscd,amssymb,latexsym}
\usepackage{color}
\usepackage{graphicx}

\newtheorem{Theorem}{\hspace{\parindent}\bf Theorem}[section]
\newtheorem{Lemma}{\hspace{\parindent}\bf Lemma}[section]
\newtheorem{Proposition}{\hspace{\parindent}\bf Proposition}[section]
\newtheorem{Corollary}{\hspace{\parindent}\bf Corollary}[section]
\newtheorem{Remark}{\hspace{\parindent}\bf Remark}[section]
\newtheorem{Definition}{\hspace{\parindent}\bf Definition}[section]
\newcommand{\pf}{\vskip0.1cm \noindent {\it Proof. }}
\newcommand{\pff}[1]{\vskip0.1cm \noindent {\it Proof of #1. }}

\def\ext{\mathop{\rm E}_\alpha\nolimits}
\def\div{\mathop{\rm div}\nolimits}

\DeclareMathOperator*{\limm}{lim}

\DeclareMathOperator*{\inff}{inf}

 \allowdisplaybreaks
\DeclareMathOperator*{\supp}{supp}
\DeclareMathOperator*{\dist}{dist}

\begin{document}

\title{\bf On some critical problems for the fractional Laplacian operator}


\begin{abstract}
We study the effect of lower order perturbations in the existence of positive
solutions to the following critical elliptic problem involving the fractional
Laplacian:
$$
   \left\{\begin{array}{ll}
    (-\Delta)^{\alpha/2}u=\lambda u^q+u^{\frac{N+\alpha}{N-\alpha}}, \quad u>0 &\quad\mbox{in }
    \Omega, \\
    u=0&\quad\mbox{on } \partial\Omega,
  \end{array}\right.
$$
where $\Omega\subset\mathbb{R}^N$ is a smooth bounded domain, $N\ge1$,
$\lambda>0$, $0<q<\frac{N+\alpha}{N-\alpha}$, $0<\alpha<\min\{N,2\}$. For
suitable conditions on $\alpha$ depending on $q$, we prove: In the case $q<1$,
there exist at least two  solutions for every $0<\lambda<\Lambda$ and some
$\Lambda>0$, at least one if $\lambda=\Lambda$, no solution if
$\lambda>\Lambda$. For $q=1$ we show existence of at least one solution for
$0<\lambda<\lambda_1$ and nonexistence for $\lambda\ge\lambda_1$. When $q>1$
the existence is shown for every $\lambda>0$. Also we prove that the solutions
are bounded and regular.
\end{abstract}

\author[B. Barrios]{B. Barrios}
\address{Departamento de Matem\'{a}ticas\\
Universidad Aut\'onoma de Madrid, 28049 Madrid, Spain. And
Instituto de Ciencias Matem\'{a}ticas, (ICMAT, CSIC-UAM-UC3M-UCM),
C/Nicol\'{a}s Cabrera 15, 28049 Madrid, Spain.}
\email{bego.barrios@uam.es}

\author[E. Colorado]{E. Colorado}
\address{Departamento de Matem{\'a}ticas\\
Universidad Carlos III de Madrid, 28911 Legan{\'e}s (Madrid),
Spain. And Instituto de Ciencias Matem\'{a}ticas, (ICMAT,
CSIC-UAM-UC3M-UCM), C/Nicol\'{a}s Cabrera 15, 28049 Madrid,
Spain.} \email{eduardo.colorado@uc3m.es}

\author[A. de Pablo]{A. de Pablo}
\address{Departamento de Matem{\'a}ticas\\
Universidad Carlos III de Madrid, 28911 Legan{\'e}s (Madrid),
Spain.} \email{arturop@math.uc3m.es}

\author[U. S{\'a}nchez]{U. S{\'a}nchez}
\address{Departamento de Matem{\'a}ticas\\
Universidad Carlos III de Madrid, 28911 Legan{\'e}s (Madrid),
Spain.} \email{urko.sanchez@gmail.com}

\thanks{B. Barrios is  partially supported by Research Projects of MICINN of
Spain (Ref. MTM2010-16518)\\
\indent E. Colorado is partially supported by Research Projects of
MEC of
Spain (Ref. MTM2009-10878, MTM2010-18128) and of Comunidad de Madrid-UC3M (Ref. CCG10-UC3M/ESP-4609).\\
\indent A. de Pablo is partially supported by Research Projects of
MEC of Spain (Ref.  MTM2008-06326-C02-02).}
\date{\today}
\maketitle

\section{Introduction}
\setcounter{equation}{0} Problems of the type
\begin{equation}\label{uno}
  \left\{\begin{array}{ll}
    -\Delta u=f(u)&\quad\mbox{in } \Omega, \\
    u=0&\quad\mbox{on } \partial\Omega,
  \end{array}\right.
\end{equation}
for different kind of nonlinearities $f$, have been the main subject of
investigation in a large amount of works in the last thirty years. See for
example the list (far from complete)~\cite{A, ABC,brezis-nirenberg,Lions85b}.
One of the most important cases of problem \eqref{uno} is the critical power
$f(u)=u^{\frac{N+2}{N-2}}$, $N>2$, since it is well known that this problem has
no positive solution  provided the domain is starshaped. In a pioneering work
\cite{brezis-nirenberg}, Brezis and Nirenberg showed that, contrary to
intuition, the critical problem with small linear perturbations can provide positive solutions. After
that, in
\cite{ABC}, using the results on concentration-compactness of
Lions, \cite{Lions85b}, Ambrosetti, Brezis and Cerami proved some
results on existence and multiplicity of solutions for a sublinear
perturbation of the critical power, among others.

Recently, several studies have been performed for classical
elliptic equations with the Laplacian operator substituted by its
fractional powers. In particular, in \cite{tan} it is studied the
problem
\begin{equation}\label{mainprob}
   \left\{\begin{array}{ll}
    (-\Delta)^{1/2}u=\lambda u+u^{\frac{N+1}{N-1}}&\quad\mbox{in } \Omega, \\
    u=0&\quad\mbox{on } \partial\Omega,
  \end{array}\right.
\end{equation}
the analogue case to the problem in \cite{brezis-nirenberg}, but with the
square root of the Laplacian instead of the Laplacian. This operator is defined
in \cite{cabre-tan} through the spectral decomposition of the Laplacian
operator in $\Omega$ with zero Dirichlet boundary conditions.  Prior to this
study, in
\cite{cabre-tan} the authors proved that there is no solution in
the case $\lambda=0$ and $\Omega$ starshaped.

In this paper we are interested in the following perturbations of
the critical power case for different powers of the Laplacian,
$$
(P_\lambda)\quad
   \left\{\begin{array}{ll}
    (-\Delta)^{\alpha/2}u=f_\lambda (u)&\quad\mbox{in } \Omega, \\
    u=0&\quad\mbox{on } \partial\Omega,
  \end{array}\right.
$$
with $f_\lambda (u):=\lambda u^q+u^{\frac{N+\alpha}{N-\alpha}}$,
$0<q<\frac{N+\alpha}{N-\alpha}$, $0<\alpha<2$ and $N>\alpha$. Along the paper
we will look only for positive solutions to $(P_\lambda)$ (so many times we will omit the term
``positive").

For the definition of the fractional Laplacian operator we follow some ideas of
\cite{cabre-tan}, together with results from
\cite{brandle-colorado-depablo-sanchez} and
\cite{caffarelli-silvestre}. In particular, we define the
eigenvalues $\lambda_j$ of $(-\Delta)^{\alpha/2}$ as the power
$\alpha/2$ of the eigenvalues $\rho_j$ of $(-\Delta)$, i.e.,
$\lambda_j=\rho_j^{\alpha/2}$; both with zero Dirichlet boundary
data. With this definition, it has been proved in
\cite{brandle-colorado-depablo-sanchez}, using a generalized
Pohozaev identity, that problem $(P_\lambda)$ has no solution for
$\lambda=0$ whenever $\Omega$ is a starshaped domain.

Our main results dealing with Problem $(P_\lambda)$ are the
following.
\begin{Theorem}\label{maintheorem2}
Let $0<q<1$. Then, there exists $0<\Lambda <\infty$ such
that the problem $(P_\lambda)$
\begin{enumerate}
\item has no positive solution  for  $\lambda>\Lambda$;
\item has a minimal positive solution for any $0<\lambda\le\Lambda$. Moreover the family of minimal solutions is
increasing with respect to $\lambda$;
\item if $\lambda=\Lambda$ there is at least one positive solution;
\item if $\alpha\ge1$ there are at least two positive solutions for $0<\lambda<\Lambda$.
\end{enumerate}
\end{Theorem}

\begin{Theorem}\label{maintheorem}
Let $q=1$,  $0<\alpha<2$ and $N\geq 2\alpha$. Then the problem
$(P_\lambda)$
\begin{enumerate}
\item has no positive solution  for  $\lambda\geq\lambda_{1}$;
\item has at least one positive solution for each
$0<\lambda<\lambda_{1}$.
\end{enumerate}
\end{Theorem}

\begin{Theorem}\label{maintheorem3}
Let $1<q<\frac{N+\alpha}{N-\alpha}$,  $0<\alpha<2$ and $N> \alpha(1+1/q)$. Then
the problem $(P_\lambda)$ has at least one positive solution  for any
$\lambda>0$.
\end{Theorem}

The restriction $\alpha\ge1$ in Theorem \ref{maintheorem2}-$(4)$ seems to be
technical. We remember that in the study of the corresponding subcritical case
performed in
\cite{brandle-colorado-depablo-sanchez} the same restriction on
$\alpha$ appeared. In that case the difficulty was to find  a
Liouville-type theorem for $0<\alpha<1$. Here, due to  the lack of
regularity, see Proposition~\ref{prop:reg}, it is not clear how to
separate the solutions in the appropriate way,
Lemma~\ref{separacion}, see also \cite{colorado,D}.

On the other hand, we have left open the range $\alpha<N<2\alpha$ in Theorem \ref{maintheorem}. See the
special case $\alpha=2$ and $N=3$ in \cite{brezis-nirenberg}. If $\alpha=1$
this range is empty, see \cite{tan}.

As to the regularity of solutions, they are bounded and
``classical" in the sense that they have as much regularity as it
is required in the equation, i.e., they possess $\alpha$
``derivatives", see Propositions~\ref{prop:bound}
and~\ref{prop:reg}. Even more, if $\alpha=1$, they belong to
$\mathcal{C}^{1,q}(\overline\Omega)$ or
$\mathcal{C}^\infty(\overline\Omega)$, whenever $0<q<1$ or
$q\ge1$, respectively.

\textsc{Organization of the paper.} In a preliminary Section~\ref{sect-preliminares} we
describe the appropriate functional setting for the study of
problem $(P_\lambda)$, including the definition of an equivalent problem,
with the aid of an extra variable, which provides some advantages, see Remark \ref{remark-diff}. Then we devote Sections~\ref{sect-sublinear}
and~\ref{sect-linear-superlinear} to the proof of
Theorems~\ref{maintheorem2}--\ref{maintheorem3}. Finally the regularity
results, together with a concentration-compactness theorem, are proved in
Section~\ref{sect-regularity}.


\section{Preliminaries and functional setting}\label{sect-preliminares}
\setcounter{equation}{0}

The powers $(-\Delta)^{\alpha/2}$ of the positive Laplace operator $(-\Delta)$,
in a bounded domain $\Omega$ with zero Dirichlet boundary data, are defined
through the spectral decomposition using the powers of the eigenvalues of the
original operator. Let $(\varphi_j,\rho_j)$ be the eigenfunctions and
eigenvectors of $(-\Delta)$ in $\Omega$ with zero Dirichlet boundary data. Then
$(\varphi_j,\rho_j^{\alpha/2})$ are the eigenfunctions and eigenvectors of
$(-\Delta)^{\alpha/2}$, also with  Dirichlet boundary conditions. In fact, the
fractional Laplacian $(-\Delta)^{\alpha/2}$ is well defined in the space of
functions
$$
H^{\alpha/2}_0(\Omega)=\left\{u=\sum a_j
  \varphi_j\in L^2(\Omega)\  :\  \|u\|_{H^{\alpha/2}_0(\Omega)}=
  \left(\sum a_j^2\rho_j^{\alpha/2}\right)^{1/2}<\infty\right\},
$$
and, as a consequence,
$$
(-\Delta)^{\alpha/2}u=\sum
  a_j\rho_j^{\alpha/2}\varphi_j\,.
$$
Note that then
$\|u\|_{H^{\alpha/2}_0(\Omega)}=\|(-\Delta)^{\alpha/4}u\|_{L^2(\Omega)}$.

The dual space $H^{-\alpha/2}(\Omega)$ is defined in the standard
way, as well as the inverse operator $(-\Delta)^{-\alpha/2}$.

We now consider the problem
\begin{equation}
\left\{\begin{array}{ll}
    (-\Delta)^{\alpha \slash 2} u=f(u) &\quad\mbox{in } \Omega, \\
    u=0&\quad\mbox{on } \partial\Omega,
  \end{array}\right.
\label{f(u)}\end{equation} in this functional framework. Since the
above definition of the fractional Laplacian allows to integrate
by parts in the proper spaces, a natural definition of energy
solution to problem \eqref{f(u)} is the following.
\begin{Definition}
We say that $u\in H^{\alpha/2}_0(\Omega)$ is a solution of
\eqref{f(u)} if the identity
\begin{equation}\label{energy_sol}
\int_{\Omega} (-\Delta)^{\alpha/4}u(-\Delta)^{\alpha/4} \varphi\,
dx = \int_{\Omega}f(u)\varphi  \, dx
\end{equation}
holds for every function $\varphi \in H^{\alpha/2}_0(\Omega)$.
\end{Definition}

Our problem $(P_\lambda)$ is like problem \eqref{f(u)} with
$f(u)=f_\lambda (u)=\lambda u^q+u^{\frac{N+\alpha}{N-\alpha}}$. In this case the
right-hand side of \eqref{energy_sol} is well defined since
$\varphi\in H^{\alpha/2}_0(\Omega)\hookrightarrow
L^{\frac{2N}{N-\alpha}}(\Omega)$, while $u\in
H^{\alpha/2}_0(\Omega)$ hence $f(u)\in
L^{\frac{2N}{N+\alpha}}(\Omega) \hookrightarrow
H^{-\alpha/2}(\Omega)$.

Associated to problem \eqref{f(u)} we consider the energy functional
$$
  I(u)=\frac12\int_\Omega\left|(-\Delta)^{\alpha/4}u\right|^2\,dx-
  \int_\Omega F(u)\,dx\,,
$$
where $F(u)=\int_0^u f(s)\,ds$. In our case it reads
\begin{equation}\label{funct-abajo}
  I(u)=\frac12\int_\Omega\left|(-\Delta)^{\alpha/4}u\right|^2\,dx-
  \frac{\lambda}{q+1}\int_\Omega u^{q+1}\,dx-
  \frac{N-\alpha}{2N}\int_\Omega u^{\frac{2N}{N-\alpha}}\,dx\,.
\end{equation}
This functional is well defined in $H_0^{\alpha/2}(\Omega)$, and moreover, the
critical points of $I$ correspond to solutions to $(P_\lambda)$.

We now include the main ingredients of a recently developed technique
used in order to deal with fractional powers of the Laplacian.

Motivated by the work of Caffarelli and Silvestre
\cite{caffarelli-silvestre}, several authors have considered an
equivalent definition of the operator $(-\Delta)^{\alpha/2}$ in a
bounded domain with zero Dirichlet boundary data by means of an
auxiliary variable, see
\cite{brandle-colorado-depablo-sanchez,cabre-sire,cabre-tan,capella-davila-dupaigne-sire,stinga-torrea}.

Associated to the bounded domain $\Omega$, let us consider the
cylinder
$\mathcal{C}_\Omega=\Omega\times(0,\infty)\subset\mathbb{R}^{N+1}_+$.
The points in $\mathcal{C}_\Omega$ are denoted by $(x,y)$. The
lateral boundary of the cylinder will be denoted by
$\partial_L\mathcal{C}_\Omega=\partial\Omega\times(0,\infty)$.
Now, for a function $u\in H_0^{\alpha/2}(\Omega)$, we define the
{\it $\alpha$-harmonic extension} $w=\ext(u)$ to the cylinder
$\mathcal{C}_\Omega$ as the solution to the problem
\begin{equation}\label{extension}
  \left\{\begin{array}{ll}
    \div(y^{1-\alpha}\nabla w)=0&\quad \mbox{in } \mathcal{C}_\Omega,  \\
    w=0&\quad \mbox{on } \partial_L\mathcal{C}_\Omega, \\
    w=u&\quad \mbox{on } \Omega\times\{y=0\}.
  \end{array}
  \right.
\end{equation}
The extension function belongs to the space
$$
X_0^\alpha(\mathcal{C}_\Omega)=\left\{z\in
L^{2}(\mathcal{C}_\Omega)\,:\: z=0 \mbox{ on }
\partial_L\mathcal{C}_\Omega,\  \|z\|_{X_0^\alpha(\mathcal{C}_{\Omega})}=
\left(\kappa_\alpha\int_{\mathcal{C}_\Omega} y^{1-\alpha}|\nabla
z|^2\right)^{1/2}<\infty\right\},
$$
where $\kappa_\alpha$ is a normalization constant. With this
constant we have that the extension operator is an isometry
between $H_0^{\alpha/2}(\Omega)$ and
$X_0^\alpha(\mathcal{C}_\Omega)$. That is
\begin{equation}\label{equivalencia_normas}
\|\ext(\psi)\|_{X^{\alpha}_{0}(\mathcal{C}_{\Omega})}=\|\psi\|_{H^{\alpha/2}_{0}(\Omega)}\,
, \quad \forall \, \psi\in H_{0}^{\alpha/2}(\Omega).
\end{equation}
Moreover, for any function $\varphi\in
X_{0}^{\alpha}(\mathcal{C}_{\Omega})$, we have the following trace
inequality
\begin{equation}\label{trace}
\|\varphi(\cdot,0)\|_{H^{\alpha/2}_{0}(\Omega)}\leq\|\varphi\|_{X^{\alpha}_{0}(\mathcal{C}_{\Omega})}.
\end{equation}
The relevance of the extension function $w$ is that it is related to
the fractional Laplacian of the original function $u$ through the
formula
\begin{equation}
    \label{normalder}
-\lim\limits_{y\to0^+}y^{1-\alpha}\dfrac{\partial w}{\partial
y}(x,y)=\frac{1}{\kappa_\alpha}(-\Delta)^{\alpha/2}u(x),
\end{equation}
see
\cite{brandle-colorado-depablo-sanchez,cabre-sire,cabre-tan,caffarelli-silvestre,capella-davila-dupaigne-sire,stinga-torrea}.
When $\Omega=\mathbb{R}^N$, the above Dirichlet-Neumann procedure
\eqref{extension}--\eqref{normalder} provides a formula for the
fractional Laplacian in the whole space equivalent to that obtained
from Fourier Transform, see \cite{caffarelli-silvestre}. In that
case, the $\alpha$-harmonic extension  and the fractional Laplacian
have explicit expressions in terms of the Poisson and the
Riesz kernels, respectively:
\begin{equation}
  \label{poisson}
\begin{array}{l}
w(x,y)=P_y^\alpha*u(x)=
c_{N,\alpha}y^{\alpha}\displaystyle\int_{\mathbb{R}^N}\dfrac{u(s)}{(|x-s|^{2}+y^{2})^{\frac{N+\alpha}{2}}}\,ds\,,
\\ [5mm]
(-\Delta)^{\alpha/2}u(x)=
d_{N,\alpha}\displaystyle P.V.\int_{\mathbb{R}^N}\frac{u(x)-u(s)}{|x-s|^{N+\alpha}}\,ds\,.
\end{array}
\end{equation}
In fact the extension technique is developed originally for the fractional
Laplacian defined in the whole space, \cite{caffarelli-silvestre}, and the
corresponding functional spaces are well defined on the homogeneous fractional Sobolev space
$\dot H^{\alpha/2}(\mathbb{R}^N)$ and the weighted Sobolev space
$X^\alpha(\mathbb{R}^{N+1}_+)$. The constants in \eqref{poisson} and
\eqref{normalder} satisfy the identity $\alpha c_{N,\alpha}\kappa_{\alpha}=d_{N,\alpha}$. Their explicit
value can be consulted for instance in
\cite{brandle-colorado-depablo-sanchez}. We will use the following notation,
$$
L_\alpha w:=-\div(y^{1-\alpha}\nabla w), \qquad \frac{\partial
w}{\partial \nu^\alpha}:=
  -\kappa_{\alpha}\lim_{y\to0^+}y^{1-\alpha}\frac{\partial w}{\partial y}\,.
$$
 With this extension, we can reformulate our problem
$(P_\lambda)$ as
$$
 (\overline P_\lambda)\quad\left\{
 \begin{array}{ll}
    L_{\alpha}w=0 &\quad \mbox{in } \mathcal{C}_{\Omega} \\ [1mm]
    w=0 & \quad \mbox{on } \partial_{L} \mathcal{C}_{\Omega}  \\
   \dfrac{\partial w}{\partial \nu^\alpha}=\lambda
w^q+w^{\frac{N+\alpha}{N-\alpha}} &\quad \mbox{in }
\Omega\times\{y=0\} .
  \end{array}
  \right.
$$
An energy solution to this problem is a function $w\in
X_0^\alpha(\mathcal{C}_\Omega)$ such that
\begin{equation}
\kappa_\alpha\int_{\mathcal{C}_\Omega} y^{1-\alpha}\langle\nabla
w,\nabla \varphi\rangle\, dxdy =\int_{\Omega}\left(\lambda
w^q+w^{\frac{N+\alpha}{N-\alpha}}\right)\varphi  \, dx,\qquad
\forall\;\varphi \in X^{\alpha}_0(\mathcal{C}_\Omega).
\label{debil}\end{equation}

For any energy solution $w\in X_0^\alpha(\mathcal{C}_\Omega)$ to
this problem, the function $u=w(\cdot,0)$, defined in the sense of
traces, belongs to the space $H_0^{\alpha/2}(\Omega)$ and is an
energy solution to problem $(P_\lambda)$. The converse is also
true. Therefore, both formulations are equivalent.

The associated energy functional to the problem $(\overline
P_\lambda)$ is
\begin{equation}\label{funct-arriba}
  J(w)=\frac{\kappa_\alpha}2\int_{\mathcal{C}_\Omega} y^{1-\alpha}|\nabla w|^2\,dxdy-
  \frac{\lambda}{q+1}\int_\Omega w^{q+1}\,dx-
  \frac{N-\alpha}{2N}\int_\Omega w^{\frac{2N}{N-\alpha}}\,dx\,.
\end{equation}

Clearly, critical points of $J$ in
$X_0^\alpha(\mathcal{C}_\Omega)$ correspond to critical points of
$I$ in $H_0^{\alpha/2}(\Omega)$. Even more, minima of $J$ also
correspond to minima of $I$, see Section~\ref{sect-sublinear}.

\begin{Remark}\label{remark-diff}
In the sequel, and in view of the above equivalence, we will use both
formulations of the problem, in $\Omega$ or in $\mathcal{C}_\Omega$, whenever
we may take some advantage. In particular, we will use the extension version
$(\overline P_\lambda)$ when dealing with the fractional operator acting on
products of functions, since it is not clear how to calculate this action. This
difficulty appears in the proof of the concentration-compactness result,
Theorem~\ref{lions}, among others.
\end{Remark}

Another tool which is very useful in what follows is the trace inequality
\begin{equation}
  \label{eq:trace-r}
  \int_{\mathcal{C}_\Omega} y^{1-\alpha}|\nabla z(x,y)|^2\,dxdy \ge
  C\left(\int_{\Omega} |z(x,0)|^{r}\,dx\right)^{2/r},
\end{equation}
for any $1\le r\le
\frac{2N}{N-\alpha}$, $N>\alpha$,  and any $z\in X_0^\alpha(\mathcal{C}_\Omega)$, where $C=C(\alpha,r,N,\Omega)>0$.
In fact it is equivalent to the fractional Sobolev inequality
\begin{equation}\label{sobolev}
 \int_{\Omega}|(-\Delta)^{\alpha/4}v|^2\,dx\ge C\left(\int_{\Omega} |v|^{r}\,dx\right)^{2/r}
\end{equation} for any $1\le r\le
\frac{2N}{N-\alpha}$, $N>\alpha$, and every $v\in H^{\alpha/2}_0(\Omega)$.
In the following we will denote the critical fractional Sobolev exponent
$2^*_\alpha=\frac{2N}{N-\alpha}$.

\noindent \begin{Remark} When
$r=2^*_\alpha$, the best constant in \eqref{eq:trace-r} will be denoted by
$S(\alpha,N)$. This constant is explicit and independent of the domain; its
exact value is
$$
\displaystyle
S(\alpha,N)=\frac{2\pi^{\frac{\alpha}{2}}\Gamma(\frac{N+\alpha}{2})\Gamma(\frac{2-\alpha}{2})
(\Gamma(\frac{N}{2}))^{\frac{\alpha}{N}}}{\Gamma(\frac{\alpha}{2})\Gamma(\frac{N-\alpha}{2})
(\Gamma(N))^{\frac{\alpha}{2}}}.
$$
It is not achieved in any bounded domain,  so we have
\begin{equation}
  \label{eq:trace-p}
  \int_{\mathbb{R}^{N+1}_+} y^{1-\alpha}|\nabla z(x,y)|^2\,dxdy \ge
  S(\alpha,N)\left(\int_{\mathbb{R}^{N}} |z(x,0)|^{\frac{2N}{N-\alpha}}\,dx\right)^{\frac{N-\alpha}N},
  \quad \forall\,z\in X^\alpha(\mathbb{R}^{N+1}_+)\,,
\end{equation}
though it is indeed achieved in that case $\Omega=\mathbb{R}^{N+1}_+$ when
$u=z(\cdot,0)$ takes the form
\begin{equation}\label{minimizers0}
u(x)=u_{\varepsilon}(x)= \frac{\varepsilon^{(N-\alpha)\slash
2}}{(|x|^2+\varepsilon^2)^{(N-\alpha)\slash
    2}},
\end{equation}
with $\varepsilon>0$ arbitrary and $z=\ext(u)$. See
\cite{brandle-colorado-depablo-sanchez} for more details. This
will be used in Sections \ref{sect-sublinear} and
\ref{sect-linear-superlinear}. The best constant in \eqref{sobolev}
when $\Omega=\mathbb{R}^{N}$ is then $\kappa_\alpha S(\alpha,N)$.
\end{Remark}
\section{Sublinear case: $0<q<1$.}\label{sect-sublinear}
We prove here Theorem \ref{maintheorem2}.  As we have said in Remark \ref{remark-diff}, there are some
points where it is difficult to work directly with the fractional Laplacian,
due to the absence of formula for the fractional Laplacian of a product.
Therefore we consider in some occasions the extended problem
$(\overline{P}_{\lambda})$.

To begin with that problem, we prove  that local minima of the functional $I$
correspond to local minima of the extended functional $J$.
\begin{Proposition}\label{min_equivalentes}
A function $u_{0}\in H_{0}^{\alpha\slash 2}(\Omega)$ is a local minimum of $I$
if and only if $w_{0}=\ext(u_{0})\in X^{\alpha}_{0}(\mathcal{C}_{\Omega})$ is a
local minimum of $J$.
\end{Proposition}
\pf
Firstly let $u_{0}\in H_{0}^{\alpha\slash 2}(\Omega)$ be a local minimum of
$I$. Suppose, by contradiction, that $w_{0}=\ext(u_{0})$ is not a local minimum
for the extended functional $J$. Then by (\ref{equivalencia_normas}) and
\eqref{trace}, we have that, for any $\varepsilon
> 0$, there exists $w_{\varepsilon} \in
X^{\alpha}_{0}(\mathcal{C}_{\Omega})$, with $\|w_{0}-w_{\varepsilon}
\|_{X^{\alpha}_{0}(\mathcal{C}_{\Omega})} < \varepsilon$, such that
$$
I(u_{0})=J(w_{0})>J(w_{\varepsilon})\geq I({z_{\varepsilon}})
$$
where $z_{\varepsilon} = w_{\varepsilon}(\cdot, 0) \in H_{0}^{\alpha\slash
2}(\Omega)$ satisfies $\|u_{0}-z_{\varepsilon}
\|_{H_{0}^{\alpha\slash 2}(\Omega)} < \varepsilon$.

On the other hand, let $w_{0}\in X^{\alpha}_{0}(\mathcal{C}_{\Omega})$ be a
local minimum of $J$. It is clear, from the definition of the extension
operator, that $w_{0}$ is $\alpha$-harmonic. So we conclude.~\hfill$\Box$

\

We return now to the original problem $(P_\lambda)$, posed at the
bottom $\Omega\times\{y=0\}$.

\begin{Lemma}\label{lem:>}
Let $\Lambda$ be defined by
$$
\Lambda =\sup\left\{ \lambda >0\ :\ \mbox{Problem $(P_\lambda)$ has
solution}\right\}.
$$
Then $0<\Lambda <\infty$.
\end{Lemma}

\pf Let $(\lambda_1,\varphi_1)$ be the first eigenvalue and a
corresponding positive eigenfunction of the fractional Laplacian in $\Omega$.
Then, using $\varphi_1$ as a test function in $(P_\lambda)$, we have that
\begin{equation}
  \label{eq:lem.>}
  \int_\Omega \left(\lambda u^q+u^{\frac{N+\alpha}{N-\alpha}}\right)\varphi_1\,dx=
  \lambda_1\int_\Omega u\varphi_1\, dx.
\end{equation}
Since there exist positive constants $c,\delta$ such that $\lambda
t^q+t^{\frac{N+\alpha}{N-\alpha}}>c\lambda^\delta t$, for any
$t>0$ we obtain from~\eqref{eq:lem.>}  that
$c\lambda^\delta<\lambda_1$ which implies $\Lambda<\infty$.

To prove $\Lambda>0$ we use the sub- and supersolution technique to construct a
solution for any small $\lambda$, see \cite{GP, ABC}. In fact a subsolution is
obtained as a small multiple of $\varphi_1$. A supersolution is a large
multiple of the function $g$ solution to
$$
 \left\{\begin{array}{ll}
    (-\Delta)^{\alpha \slash 2} g=1 &\quad\mbox{in } \Omega, \\
    g=0&\quad\mbox{on } \partial\Omega.
  \end{array}\right.
$$~\hfill$\Box$

Comparison is clear for linear problems associated to the
fractional Laplacian, as it is for the Laplacian. On the other
hand, it is in general not true for nonlinear problems.
Nevertheless, it holds when the reaction term is a nonnegative
sublinear function, see
\cite{BK,ABC,brandle-colorado-depablo-sanchez}. Therefore, it is easy
to show, comparing with the problem with only the concave terms
$\lambda u^q$, that in fact there is at least one positive
solution $u_\lambda$ to problem $(P_\lambda)$ for every $\lambda$
in the whole interval $(0,\Lambda)$. Even more, these constructed
solutions are minimal and are increasing with respect to $\lambda$
(see Lema 5.7 of \cite{brandle-colorado-depablo-sanchez}).

To prove existence of solution in the extremal value $\lambda=\Lambda$, the
idea, like in
\cite{ABC}, consists on passing to the limit as $\lambda_n\nearrow\Lambda$ on
the sequence $\{z_n\}=\{z_{\lambda_{n}}\}$, where $z_{\lambda_{n}}$ is the
minimal solution of $(\overline{P}_{\lambda})$ with $\lambda=\lambda_{n}$.
Denote by $J_{\lambda_n}$ the associated functional. Clearly
$J_{\lambda_{n}}(z_n)< 0$, hence
$$
 \displaystyle 0> J_{\lambda_{n}}(z_n) - \frac{1}{2^{*}_{\alpha}}\langle J'_{\lambda_{n}}(z_n),z_n\rangle
 =
\kappa_\alpha\left(\frac 12-\frac{1}{2^*_\alpha}\right)\|  z_n\|^{2}_{X^{\alpha}_{0}(\mathcal{C}_{\Omega})}
-\lambda_n\left(\frac{1}{q+1}-\frac{1}{2^*_\alpha}\right) \int_\Omega
z_n^{q+1}dx.
$$
Therefore, by the Sobolev  and Trace inequalities, \eqref{sobolev} and
\eqref{trace} respectively, there exits a constant $C>0$ such that $\|
z_n\|_{X_0^\alpha(\mathcal{C}_{\Omega})}\le C.$ As a consequence, there exists
a subsequence weakly convergent to some $z_\Lambda$ in
$X_0^\alpha(\mathcal{C}_{\Omega})$. By comparison, $z_\Lambda\ge z_\lambda>0$,
for any $0<\lambda<\Lambda$, so one gets easily that $z_\Lambda$ is a weak
nontrivial solution to $(\overline{P}_{\lambda})$ with $\lambda=\Lambda$.

Having proved the first three items in Theorem~\ref{maintheorem2}, we focus in
the sequel on proving the existence of a second solution, for which we recall that $\alpha\ge 1$.

The proof is divided into several steps: we first show that the minimal solution is a local minimum for the functional
 $I$; so we can use the Mountain Pass Theorem, obtaining a minimax Palais-Smale (PS) sequence.
 In the next step, in order to find a second solution, we prove a local (PS)$_c$ condition for $c$ under
 a critical level $c^*$. To do that,  we will construct
path by localizing the minimizers of the Trace/Sobolev inequalities at the possible Dirac
Deltas, given by the concentration-compactness result in Theorem \ref{lions}.

We begin with a separation lemma in the $\mathcal{C}^1$-topology.
\begin{Lemma}\label{separacion}
Let $0<\mu_1<\lambda_0<\mu_2<\Lambda$. Let $z_{\mu_1}$,
$z_{\lambda_0}$ and $z_{\mu_2}$ be the corresponding minimal
solutions to $(P_{\lambda})$, $\lambda=\mu_1,\,\lambda_0$ and
$\mu_2$ respectively. If $X=\{z\in \mathcal{C}_0^1(\Omega)|
\,z_{\mu_{1}}\leq z \leq z_{\mu_{2}} \}$, then there exists
$\varepsilon
>0$ such that
$$
\{z_{\lambda_0}\} + \varepsilon B_{1} \subset X,
$$
where $B_{1}$ is the unit ball in $\mathcal{C}_0^1(\Omega)$.
\end{Lemma}
\pf Since $\alpha\ge1$, we have that any solution $u$ to
$(P_{\lambda})$, for arbitrary $0<\lambda<\Lambda$ belongs to
$\mathcal{C}^{1,\gamma}(\overline{\Omega})$ for some positive $\gamma$, see
Proposition~\ref{prop:reg}. Therefore, we deduce that there exists
a positive constant $C$ such that
\begin{equation}\label{distancia2}
u(x)\leq C\dist(x,\partial\Omega),\,x\in\Omega.
\end{equation}
On the other hand, by comparison with the first eigenfunction of
the fractional Laplacian (which is indeed the first eigenfunction
$\varphi_{1}$ of the classical Laplacian), we get that there
exists a positive constant $c$ such that
\begin{equation}\label{distancia1}
u(x)\geq c\dist(x,\partial\Omega),\,x\in\Omega.
\end{equation}
These two estimates jointly with the regularity implies the result of  the
lemma. ~\hfill$\Box$

\

With this result we now obtain  a local minimum of the functional
$I$ in $\mathcal{C}_{0}^{1}(\Omega)$, as a first step, to obtain a
local minimum in $H_{0}^{\alpha/2}(\Omega)$.
\begin{Lemma}\label{minimum}
For all $\lambda \in (0,\Lambda)$ there exists a solution for
$(P_{\lambda})$ which is a local minimum of the functional $I$ in
the $\mathcal{C}^1$-topology.
\end{Lemma}
\pf Given $0<\mu_{1}<\lambda<\mu_{2}<\Lambda$, let $z_{\mu_{1}}$
and $z_{\mu_{2}}$ be the minimal solutions of $(P_{\mu_{1}})$ and
$(P_{\mu_{2}})$ respectively. Let $z:=z_{\mu_{2}}-z_{\mu_{1}}$.
Since $z_{\mu_{1}}$ and $z_{\mu_{2}}$ are properly ordered, then
$$
\left\{
 \begin{array}{ll}
    (-\Delta)^{\alpha/2}z\geq0&\quad \mbox{in } \Omega, \\
    z=0 & \quad \mbox{on } \partial \Omega.
  \end{array}
  \right.
$$
We set
$$
f^*(x,s)= \left\{
\begin{array}{ll}
f_{\lambda}(z_{\mu_{1}}(x)) &\mbox{ if } s\leq z_{\mu_{1}} ,\\ [2mm]
f_{\lambda}(s) &\mbox{ if } z_{\mu_{1}}\leq s \leq z_{\mu_{2}},\\ [2mm]
f_{\lambda}(z_{\mu_{2}}(x)) &\mbox{ if } z_{\mu_{2}}\le s,
\end{array}\right.
$$
$$
F^*(x,z)= \int_{0}^{z}f^*(x,s)\,ds
$$
and
$$
I^*(z)= \frac{1}{2}\|z\|_{H_{0}^{\alpha\slash 2}(\Omega)} -
\int_{\Omega}F^*(x,u)dx.
$$
Standard calculation shows that $I^*$ achieves its global minimum
at some $u_{0} \in H_{0}^{\alpha\slash 2}(\Omega)$, that is
\begin{equation}\label{diferencia1}
I^*(u_{0})\leq I^*(z)\quad\forall\, z \in H_{0}^{\alpha\slash
2}(\Omega).
\end{equation}
 Moreover it holds
$$
\left\{
 \begin{array}{ll}
    (-\Delta)^{\alpha/2}u_{0}= f^*(x,u_{0})&\quad \mbox{in } \Omega, \\
    u_{0}=0 & \quad \mbox{on } \partial \Omega.
  \end{array}
  \right.
$$
By Lemma \ref{separacion}, it follows that $\{u_{0}\}+\varepsilon
B_{1}\subseteq X$
for $0<\varepsilon$ small enough. Let now $z$ satisfying
$$\| z-u_{0}\|_{\mathcal{C}_{0}^{1}(\Omega)}\leq\frac{\varepsilon}{2}.$$
As $I^*(z)-I(z)$ is zero for every $z$ such that
$\|z-u_{0}\|_{\mathcal{C}_{0}^{1}(\Omega)}\leq\frac{\varepsilon}{2}$,
by (\ref{diferencia1}) we obtain that
\[I(z)=I^*(z)\geq I^*(u_{0})=I(u_{0}),\quad \forall\; z\in \mathcal{C}_{0}^{1}(\Omega),\,\mbox{ with }\|z-u_{0}\|_{\mathcal{C}_{0}^{1}(\Omega)}\leq\frac{\varepsilon}{2}.\]
~\hfill$\Box$

To show that we have obtained the desired minimum in $H_{0}^{
\alpha/2}(\Omega)$, we now check that the result by Brezis and
Nirenberg in \cite{brezis-nirenberg-H1} is also valid in our
context.
\begin{Proposition}\label{alpha_vs_c1}
Let $z_{0}\in H^{\alpha\slash 2}_0(\Omega)$ be a local minimum of
$I$ in $\mathcal{C}_{0}^{1}(\Omega)$, i.e., there exists $r>0$ such
that
\begin{equation}\label{hhip}
I(z_{0})\leq I(z_{0}+z) \qquad \forall z\in
\mathcal{C}_{0}^{1}(\Omega)\mbox{ with }
\left\|z\right\|_{\mathcal{C}_{0}^{1}(\Omega)}\leq r.
\end{equation}
Then $z_{0}$ is a local minimum of $I$ in $H_{0}^{\alpha\slash
2}(\Omega)$, that is, there exists $\varepsilon_{0}>0$ such that
$$
I(z_{0})\leq I(z_{0}+z) \qquad \forall z\in H_{0}^{\alpha\slash
2}(\Omega) \mbox{ with } \|z\|_{H_{0}^{\alpha\slash
2}(\Omega)}\leq \varepsilon_{0}.
$$
\end{Proposition}
\pf Arguing  by contradiction we suppose that
\[\forall\,\varepsilon>0,\,\exists\,z_{\varepsilon}\in
B_{\varepsilon}(z_{0}) \quad\mbox{such that }\quad
I(z_{\varepsilon})<I(z_{0}),\]
\noindent where $B_{\varepsilon}(z_{0})=\left\{z\in H_0^{\alpha/
2}(\Omega):\, \,\|z-z_{0}\|_{H_0^{\alpha/
2}(\Omega)}\leq\varepsilon\right\}$.

\

For every $j>0$ we consider the truncation map given by
$$
   T_{j}(r)\equiv
   \left\{\begin{array}{ll}
    r&\quad 0<r<j, \\
    j&\quad r\geq j.
 \end{array}\right.
 $$
Let
$$
f_{\lambda,j}(s)=f_{\lambda}(T_{j}(s)),\qquad
F_{j}(s)=\int_{0}^{u}{f_{\lambda,j}(s)ds}\,,\,u>0\,,
$$
and
\[I_{j}(z)=\frac{1}{2}\|z\|_{H_0^{\alpha/2}(\Omega)}^{2}-\int_{\Omega}{F_{j}(z)dx}.\]
Note that for each $z\in H^{\alpha/2}_{0}(\Omega)$ we have that
$I_{j}(z)\to I(z)$ as $j\to\infty$. Hence, for each $\varepsilon>0$
there exists $j(\varepsilon)$ big enough such that
$I_{j(\varepsilon)}(z_{\varepsilon})<I(z_{0})$. Clearly
$\displaystyle\min_{B_{\varepsilon}(z_{0})}{I_{j(\varepsilon)}}$ is
attained at some point, say $v_{\varepsilon}$. Thus we have
$$
I_{j(\varepsilon)}(v_{\varepsilon})\le
I_{j(\varepsilon)}(z_{\varepsilon})<I(z_{0}).
$$
Now we want to prove that $v_{\varepsilon}\rightarrow z_{0}$ in
$\mathcal{C}_{0}^{1}(\Omega)$ as $\varepsilon\searrow 0$. The Euler-Lagrange
equation satisfied by $v_\varepsilon$ involves a Lagrange multiplier
$\xi_\varepsilon$ in such a way that
\begin{equation}\label{eqeuler}
\langle
I_{j(\varepsilon)}'(v_{\varepsilon}),\varphi\rangle_{H^{-\alpha/
2}(\Omega),H_0^{\alpha/ 2}(\Omega)}=\xi_{\varepsilon}\langle
v_{\varepsilon},\varphi\rangle_{H_0^{\alpha/ 2}(\Omega)},\quad
\forall\,\varphi\in H_0^{\alpha/ 2}(\Omega).
\end{equation}
Since $v_{\varepsilon}$ is a minimum of $I_{j(\varepsilon)}$, it
holds
\begin{equation}\label{a-cero}
\xi_{\varepsilon}=\frac{\langle
I_{j(\varepsilon)}'(v_{\varepsilon}),
v_{\varepsilon}\rangle}{\|v_{\varepsilon}\|^{2}_{H_{0}^{\alpha\slash
2}(\Omega)}}\leq 0\quad\mbox{ for } 0<\varepsilon\ll1,
\quad\mbox{and}\quad \xi_\varepsilon\to 0 \mbox{ as }
\varepsilon\searrow 0.
\end{equation}
Note that by (\ref{eqeuler}), 
$v_{\varepsilon}$ satisfies the problem
$$
\left\{
 \begin{array}{ll}
    (-\Delta)^{\alpha/2}v_{\varepsilon}=\frac{1}{1-\xi_{\varepsilon}}f_{\lambda,j(\varepsilon)}(v_{\varepsilon}):=
    f_{\lambda,j(\varepsilon)}^{\varepsilon}(v_{\varepsilon})&\quad \mbox{in } \Omega, \\ [2mm]
    v_{\varepsilon}=0 & \quad \mbox{on } \partial \Omega.
  \end{array}
  \right.
$$
Clearly $\|v_{\varepsilon}\|_{H_0^{\alpha\slash 2}(\Omega)}\leq
C$, thus, by Proposition \ref{prop:bound}, this implies that
$\|v_{\varepsilon}\|_{L^{\infty}(\Omega) }\leq C$. Moreover, by
\eqref{a-cero} it follows that
$\|f_{\lambda,j(\varepsilon)}^{\varepsilon}(v_{\varepsilon})\|_{L^{\infty}(\Omega)}
\leq C$. Therefore, following the proof of Proposition
\ref{prop:reg}, we get that
$\|v_{\varepsilon}\|_{\mathcal{C}^{1,r}(\overline{\Omega})}\leq C$, for
$r=\min\{q,\alpha-1\}$ and $C$ independent of $\varepsilon$. By
Ascoli-Arzel\'{a} Theorem there exists a subsequence, still
denoted by $v_{\varepsilon}$, such that
$v_{\varepsilon}\rightarrow z_{0}$ uniformly in
$\mathcal{C}_{0}^{1}(\Omega)$ as $\varepsilon\searrow 0$.
This implies that for $\varepsilon$ small enough,
 $$
 I(v_{\varepsilon})=I_{j(\varepsilon)}(v_{\varepsilon})< I(z_{0})
 $$
 for any $v_{\varepsilon}$ with $\|v_{\varepsilon}-z_{0}\|_{\mathcal{C}_{0}^{1}(\Omega)}<\varepsilon$. ~\hfill$\Box$
\\Lemma \ref{minimum} and Proposition \ref{alpha_vs_c1} provide us a local minimum in $H_{0}^{\alpha/2}(\Omega)$, which will be
denoted by $u_{0}$. We now perform a traslation in order to
simplify the calculations.

We consider the functions
\begin{equation}
g(x,s)=\left\{\begin{array}{ll} \lambda
(u_{0}+s)^{q}-\lambda u_{0}^{q} +(u_{0}+s)^{2^*_{\alpha} -1}- u_{0}^{2^*_{\alpha} -1} &  \mbox{ if }s\ge 0,\\
0 &  \mbox{ if } s<0,
\end{array}
\right.
\label{g(x,s)}\end{equation}
\begin{equation}\label{G}
G(u)=\int_{0}^{u}g(x,s)\,ds,
\end{equation}
and the energy functional
\begin{equation}\label{II}
\widetilde{I}(u)= \frac{1}{2}\|u\|^{2}_{H^{\alpha/2}_{0}(\Omega)}
- \int_{\Omega}G(x,u)dx.
\end{equation}
Since $u\in H^{\alpha/2}_{0}(\Omega)$, $G$ is well defined and bounded from
below. Let the moved problem
$$
   (\widetilde{P}_{\lambda})\quad\left\{\begin{array}{ll}
    (-\Delta)^{\alpha/2}u= g(x,u)&\quad\mbox{in }
    \Omega\subset\mathbb{R}^N,\,\lambda>0 \\
    u=0&\quad\mbox{on } \partial\Omega.
  \end{array}\right.
$$
Hence, by standard variational theory, we know that if $\widetilde{u}
\not\equiv 0$ is a critical point of $\widetilde{I}$ then it is a solution of
$(\widetilde{P}_{\lambda})$ which, by the Maximum Principle (Lemma 2.3 of
\cite{capella-davila-dupaigne-sire}), it is $\widetilde{u}>0$. Therefore
$u=u_{0}+\widetilde{u}$ will be a second solution of $(P_{\lambda})$ for the
sublinear case. Thus we will need to study the existence of these non-trivial
critical points for $\tilde{I}$.
\\
Firstly we have
\begin{Lemma}\label{lemma:minimo}
$u=0$ is a local minimum of $\widetilde{I}$ in
$H^{\alpha/2}_{0}(\Omega).$
\end{Lemma}
\pf The proof follows the lines of \cite{ABC}, so we will be brief
in details. Note that by Proposition \ref{alpha_vs_c1} it is
sufficient to prove that $u=0$ is a local minimum of
$\widetilde{I}$ in $\mathcal{C}_{0}^{1}(\Omega)$.

Let $u\in \mathcal{C}_{0}^{1}(\Omega)$, then
\begin{equation}
G(u)=F(u_{0}+u)-F(u_{0})-\left(\lambda
u_{0}^{q}+u_{0}^{2^*_{\alpha}-1}\right)u.
\end{equation}
Therefore
$$
\begin{array}{rcl}
\displaystyle \widetilde{I}(u)&=&
\displaystyle\frac{1}{2}\|u^{}\|_{H^{\alpha/2}_{0}(\Omega)}^{2}-
\int_{\Omega}G(u)dx \\ [3mm]
&=&\displaystyle\frac{1}{2}\|u^{}\|_{H^{\alpha/2}_{0}(\Omega)}^{2} -
\int_{\Omega}F(u_{0} + u^{})dx + \int_{\Omega}F(u_{0})dx +
\int_{\Omega}\left(\lambda u_{0}^{q}+u_{0}^{2^*_{\alpha}-1}\right)u^{}dx.
\end{array}
$$
On the other hand,
$$
\begin{array}{rcl}
\displaystyle I(u_{0} + u^{})&=& \displaystyle\frac{1}{2}\|u_{0} +
u^{}\|_{H^{\alpha/2}_{0}(\Omega)}^{2} -
\int_{\Omega}F(u_{0} + u^{})dx \\ [3mm]
&=&\displaystyle\frac{1}{2}\|u_{0}\|_{H^{\alpha/2}_{0}(\Omega)}^{2} +
\frac{1}{2}\|u^{}\|_{H^{\alpha/2}_{0}(\Omega)}^{2} +
\int_{\Omega}(-\Delta)^{\alpha/4}u_{0}(-\Delta)^{\alpha/4}u^{}dx -
\int_{\Omega}F(u_{0} + u^{})dx\\ [3mm]
&=&\displaystyle\frac{1}{2}\|u_{0}\|_{H^{\alpha/2}_{0}(\Omega)}^{2} +
\frac{1}{2}\|u^{}\|_{H^{\alpha/2}_{0}(\Omega)}^{2} +
\int_{\Omega}\left(\lambda u_{0}^{q}+u_{0}^{2^*_{\alpha}-1}\right)u^{}dx -
\int_{\Omega}F(u_{0} + u^{})dx.
\end{array}
$$
Finally, as $u_{0}$ is a local minimum of $I$, we have that
$$
\begin{array}{rcl}
\displaystyle \widetilde{I}(u)&=&\displaystyle I(u_{0} + u^{}) -
\frac{1}{2}\|u_{0}\|_{H^{\alpha/2}_{0}(\Omega)}^{2} +
\int_{\Omega}F(u_{0})dx \\ [3mm]
&=&\displaystyle I(u_{0} + u^{}) - I(u_{0})\\ [3mm]&\geq& 0=\widetilde{I}(0)
\end{array}
$$
provided $\|u\|_{\mathcal{C}_{0}^{1}(\Omega)}<\varepsilon.$
~\hfill$\Box$

As a consequence of Proposition~\ref{min_equivalentes}, we obtain for the moved
functional
$$
\widetilde{J}(w)=\frac{1}{2}\|w\|^{2}_{X^{\alpha}_{0}(\mathcal{C}_{\Omega})}-\int_{\Omega}G(w(x,0))dx,
$$
with $G$ as in \eqref{g(x,s)}-\eqref{G}, the following result.
\begin{Corollary}\label{corol}
$w=0$ is a local minimum of $\widetilde{J}$ in
$X^{\alpha}_{0}(\mathcal{C}_{\Omega}).$
\end{Corollary}

Now assuming that $v=0$ is
the unique critical point of the moved functional $\widetilde{J}$, then a local
 (PS)$_{c}$ condition can be proved for $c$ under a critical level
$c^*$,
\begin{equation}\label{critical c}
c^*=\frac{\alpha}{2N}(\kappa_{\alpha}S(\alpha,N))^{\frac{N}{\alpha} }.
\end{equation}
Following the ideas given in \cite{ABC}, and by an extension
 of a concentration-compactness result by Lions,
that we prove in Theorem \ref{lions}, we obtain the following result.
\begin{Lemma}\label{strongly_convergent} If $v=0$ is the only
critical point of $\widetilde{J}$ in
$X^{\alpha}_{0}(\mathcal{C}_{\Omega})$ then $\widetilde{J}$
satisfies a local Palais Smale condition below the critical level
$\textcolor{black}{c^*}$.
\end{Lemma}
\pf Let $\{w_{n}\}$ be a Palais-Smale sequence for $\widetilde{J}$
verifying
\begin{equation}\label{cc}
\widetilde{J}(w_{n})\to
c<\textcolor{black}{c^*},\qquad\widetilde{J}'(w_{n})\to 0.
\end{equation}
Since the fact that $w_0$ is a critical point implies
$\widetilde{J}(w_{n})=J(z_n)-J(w_0)$, where $z_{n}=w_{n}+w_{0}$, we have that
\begin{equation}
J(z_{n})\to c+J(w_{0}),\qquad J'(z_{n})\to 0.\label{PS_cond}
\end{equation}
On the other hand, from  \eqref{cc} we get that the sequence
$\{z_n\}$ is uniformly bounded in $X^{\alpha}_{0}(\mathcal{C}_{\Omega})$. As a
consequence, up to a subsequence,
\begin{eqnarray}
\displaystyle z_{n}&\rightharpoonup& z \qquad\qquad
\mbox{ weakly in } X^{\alpha}_{0}(\mathcal{C}_{\Omega}) \nonumber\\
\displaystyle z_{n}(\cdot,0)&\to& z(\cdot,0) \qquad
\mbox{ strong in } L^{r}(\Omega), \quad \forall \, 1\leq r < 2^*_{\alpha} \label{estrella}\\
\displaystyle z_{n}(\cdot,0)&\to& z(\cdot,0) \qquad
\mbox{ a.e. in
} \Omega.\nonumber
\end{eqnarray}

 Note that as $v=0$ is the unique critical point of $\widetilde{J}$
 then, $z=w_{0}$.

In order to apply the concentration-compactness result, Theorem
\ref{lions}, first we prove the following.
\begin{Lemma}\label{lem:tight}
The sequence $\left\{y^{1-\alpha}|\nabla z_{n}|^{2}\right\}_{n\in\mathbb{N}}$ is tight, i.e., for any $\eta>0$ there exists $\rho_0>0$ such that
\begin{equation}\label{h}
\int_{\{y>\rho_0\}}{\int_{\Omega}{y^{1-\alpha}|\nabla z_{n}|^{2}dxdx}}\le \eta,\quad \forall\, n\in\mathbb{N}.
\end{equation}
\end{Lemma}
\pf
The proof of this lemma follows some arguments of Lema 2.2 in \cite{aap}. By contradiction, we suppose that there exits $\eta_{0}>0$ such that, for any $\rho>0$ one has, up to a subsequence,
\begin{equation}\label{hhh}
\int_{\{y>\rho\}}{\int_{\Omega}{y^{1-\alpha}|\nabla z_{n}|^{2}dxdy}}>\eta_{0}\,\quad\mbox{for every }\, n\in\mathbb{N}.
\end{equation}
Let $\varepsilon>0$ be fixed (to be precised later), and let $r>0$ be  such that
\[\int_{\{y>r\}}{\int_{\Omega}{y^{1-\alpha}|\nabla z|^{2}dxdy}}<\varepsilon.\]
Let $j=\left[\frac{M}{\kappa_{\alpha}\varepsilon}\right]$ be the integer part
and $I_{k}=\{y\in\mathbb{R}^{+}:\, \,r+k\leq y\leq r+k+1\}$,
$k=0,\,1,\,\ldots,\,j$. Since
$\|z_{n}\|_{X^{\alpha}_{0}(\mathcal{C}_{\Omega})}\leq M$, we clearly obtain
that
$$\sum_{k=0}^{j}{\int_{I_{k}}{\int_{\Omega}{y^{1-\alpha}|\nabla z_{n}|^{2}dxdy}}}\leq \int_{\mathcal{C}_\Omega}y^{1-\alpha}|\nabla z_{n}|^{2}dxdy\leq\varepsilon(j+1).$$
Therefore there exists $k_{0}\in\{0,\,\ldots,\,j\}$ such that (again up to a
subsequence)
\begin{equation}\label{tight}
\int_{I_{k_{0}}}{\int_{\Omega}{y^{1-\alpha}|\nabla z_{n}|^{2}dxdy}}\leq\varepsilon,\quad \forall\, n.
\end{equation}
Let $\chi\ge 0$ be the following regular non-decreasing cut-off function
  $$
  \chi(y)=\left\{
 \begin{array}{ll}
    0&\quad\mbox{if } y\leq r+k_{0}, \\
    1&\quad\mbox{if } y>r+k_{0}+1,
\end{array}\right.
$$
Define
$v_{n}(x,y)=\chi(y)z_{n}(x,y)$. Since $v_{n}(x,0)=0$  it follows that
\begin{eqnarray*}
|\langle J'(z_{n})-J'(v_{n}),v_{n}\rangle|&=&\kappa_{\alpha}\int_{\mathcal{C}_{\Omega}}
{y^{1-\alpha}\langle \nabla (z_{n}-v_{n}),\nabla v_{n}\rangle dxdy}\\
&=&\kappa_{\alpha}\int_{I_{k_{0}}}{\int_{\Omega}{y^{1-\alpha}\langle \nabla (z_{n}-v_{n}),
\nabla v_{n}\rangle dxdy}}.
\end{eqnarray*}
Moreover by the Cauchy-Schwartz inequality, (\ref{tight}) and the compact inclusion
$H^{1}(I_{k_{0}}\times\Omega,y^{1-\alpha})$ into $L^{2}(I_{k_{0}}
\times\Omega, y^{1-\alpha})$, we have
\begin{eqnarray*}
|\langle J'(z_{n})-J'(v_{n}),v_{n}\rangle| &\leq &\kappa_{\alpha}\left(\int_{I_{k_{0}}}{\int_{\Omega}{y^{1-\alpha}
|\nabla(z_{n}-v_{n})|^{2}dxdy}}\right)^{\!\frac 12}\!\left(\int_{I_{k_{0}}}{\int_{\Omega}{y^{1-\alpha}|\nabla v_{n}|^{2}dxdy}}\right)^{\!\frac12}\\
&\leq& C\,\kappa_{\alpha}\,\varepsilon.
\end{eqnarray*}
On the other hand, by (\ref{PS_cond}), we get
\[|\langle J'(v_{n}),v_{n}\rangle|\leq C\,\kappa_{\alpha}\,\varepsilon+o(1).\]
So, for $n$ sufficiently large,
\[\int_{\{y>r+k_{0}+1\}}{\int_{\Omega}{y^{1-\alpha}|\nabla z_{n}|^{2}dxdy}}
\leq\int_{\mathcal{C}_{\Omega}}{y^{1-\alpha}|\nabla v_{n}|^{2}dxdy}=
\frac{\langle J'(v_{n}),v_{n}\rangle}{\kappa_{\alpha}}\leq C\,\varepsilon.\]
This is a contradiction with \eqref{hhh}, which proves Lemma \ref{lem:tight}.~\hfill$\Box$

\pff{Lemma \ref{strongly_convergent} (cont.)}
In view of the previous result we can apply  Theorem
\ref{lions}. Therefore,  up to a subsequence, there exists an index set $I$, at most countable,
a sequence of points $\{ x_k\}\subset\Omega$, and nonnegative real numbers
$\mu_k,\,\nu_k$, such that
\begin{equation}\label{PS1}
y^{1-\alpha}|\nabla z_{n}|^{2} \rightarrow \mu \geq
y^{1-\alpha}|\nabla w_{0}|^{2} + \sum_{k\in I}
\mu_{k}\delta_{x_{k}}
\end{equation}
and
\begin{equation}\label{PS2}
|z_{n}(\cdot,0)|^{2^*_{\alpha}} \to
\nu=|w_{0}(\cdot,0)|^{2^*_{\alpha}} + \sum_{k\in I}
\nu_{k}\delta_{x_{k}}
\end{equation}
in the sense of measures, satisfying also the  relation $
\mu_{k}\geq S(\alpha,N)\nu_{k}^{\frac{2}{2^*_{\alpha}}},\quad\mbox{for every}\, k\in
I$.

We fix any $k_0\in I$, and let $\phi\in
\mathcal{C}_0^{\infty}(\mathbb{R}^{N+1}_{+})$ be a nonincreasing cut-off
function verifying $
\phi = 1$ in  $B_{1}^{+}(x_{k_{0}})$, $
\phi = 0$ in  $B_{2}^{+}(x_{k_{0}})^{c}$. Let now
$\phi_{\varepsilon}(x,y)=\phi(x/\varepsilon, y/\varepsilon)$, clearly
$|\nabla\phi_{\varepsilon}| \leq \frac{C}{\varepsilon}$. We denote
$\Gamma_{2\varepsilon}=B_{2\varepsilon}^{+}(x_{k_{0}})\cap\{y=0\}$. Then, using
$\phi_{\varepsilon}z_{n}$ as a test function in (\ref{PS_cond}), we have
\begin{eqnarray*}
& &\kappa_{\alpha}\limm_{n\to\infty}}\int_{\mathcal{C}_{\Omega}}{y^{1-\alpha}
\langle\nabla z_{n}, \nabla
\phi_{\varepsilon}\rangle z_{n} dxdy\\
&=&\limm_{n\to\infty}\left(\int_{\Gamma_{2\varepsilon}}{|z_{n}|^{2^*_{\alpha}}
\phi_{\varepsilon} \,
dx}+\lambda\int_{\Gamma_{2\varepsilon}}{|z_{n}|^{q+1}
\phi_{\varepsilon} \,
dx}-\textcolor{black}{\kappa_{\alpha}}\int_{B_{2\varepsilon}^{+}(x_{k_{0}})}{y^{1-\alpha}|\nabla
z_{n}|^{2} \phi_{\varepsilon} \, dxdy}\right).
\end{eqnarray*}
By (\ref{estrella}), (\ref{PS1}) and (\ref{PS2}) we get
\begin{equation}\label{lions1}
\begin{array}{ll}
&\displaystyle \lim_{n\to\infty}\textcolor{black}{\kappa_{\alpha}}\int_{\mathcal{C}_{\Omega}}y^{1-\alpha}
\langle\nabla z_{n}, \nabla \phi_{\varepsilon}\rangle z_{n} \,
dxdy\\
=&\displaystyle\int_{\Gamma_{2\varepsilon}}\phi_{\varepsilon} \,
d\nu+\lambda\int_{\Gamma_{2\varepsilon}}|w_{0}|^{q+1}
\phi_{\varepsilon} \,
dx-\textcolor{black}{\kappa_{\alpha}}\int_{B_{2\varepsilon}^{+}(x_{k_{0}})}\phi_{\varepsilon}
\, d\mu.
\end{array}
\end{equation}
On the other hand, using Theorem 1.6
in \cite{Fabes-Kenig-Serapioni}, with $w=y^{1-\alpha}\in A_{2}$
and $k=1$, we obtain that
\begin{eqnarray*}
\left(\int_{B_{2\varepsilon}^{+}(x_{k_{0}})}{y^{1-\alpha}|\nabla
\phi_{\varepsilon}|^{2}|z_{n}|^{2}dxdy}\right)^{1/2}&\leq&\frac{2}{\varepsilon}\left(\int_{B_{2\varepsilon}^{+}(x_{k_{0}})}{y^{1-\alpha}|z_{n}|^{2}dxdy}\right)^{1/2}\\
&\leq& C\left(\int_{B_{2\varepsilon}^{+}(x_{k_{0}})}{y^{1-\alpha}|\nabla
z_{n}|^{2}dxdy}\right)^{1/2}.
\end{eqnarray*}
Since $z_{n}\in X^{\alpha}_{0}(\mathcal{C}_{\Omega})$, the last
expression goes to zero as $\varepsilon\rightarrow 0$. Therefore
\begin{eqnarray*}
0&\leq&\limm_{n\to\infty} \left|\int_{\mathcal{C}_{\Omega}}{y^{1-\alpha}{\langle\nabla z_{n}, \nabla \phi_{\varepsilon} \rangle} z_{n} dxdy} \right|\\
&\leq &\limm_{n\to
\infty}\left(\int_{\mathcal{C}_{\Omega}}{y^{1-\alpha}|\nabla
z_{n}|^{2}dxdy}\right)^{1/2}
\left(\int_{B_{2\varepsilon}^{+}(x_{k_{0}})}{y^{1-\alpha}|\nabla
\phi_{\varepsilon}|^{2}|z_{n}|^{2}dxdy}\right)^{1/2}\longrightarrow0.
\end{eqnarray*}
Hence, by (\ref{lions1}), it follows that
$$\lim_{\varepsilon \to 0}
\left[\int_{\Gamma_{2\varepsilon}}\phi_{\varepsilon} \,
d\nu+\lambda\int_{\Gamma_{2\varepsilon}}|w_{0}|^{q+1}
\phi_{\varepsilon} \,
dx-\textcolor{black}{\kappa_{\alpha}}\int_{B_{2\varepsilon}^{+}(x_{k_{0}})}\phi_{\varepsilon}
\, d\mu \right]
=\nu_{k_{0}}-\textcolor{black}{\kappa_{\alpha}}\mu_{k_{0}}=0.$$
Therefore we get
that
$$
\nu_{k_0}=0 \qquad \mbox{ or } \qquad \nu_{k_0}\geq
(\textcolor{black}{\kappa_{\alpha}}S(\alpha,N))^{\frac{N}{\alpha}}.
$$
Suppose that $\nu_{k_{0}}\neq 0$. It follows that
\begin{eqnarray*}
c+J(w_{0})&=&\limm_{n\to\infty}{J(z_{n})- \frac{1}{2}\langle
J'(z_{n}),z_{n}\rangle}\\
&\geq& \frac{\alpha}{2N}\int_{\Omega}{w_{0}^{2^*_{\alpha}}dx} +\frac{\alpha}{2N}\nu_{k_0}+\lambda
\left(\frac{1}{2}-\frac{1}{q+1}
 \right)\int_{\Omega}{w_{0}^{q+1}dx}\\
 &\ge &J(w_{0})+\frac{\alpha}{2N}(\textcolor{black}{\kappa_{\alpha}}S(\alpha,N))^{N/\alpha} =J(w_0)+c^*.
\end{eqnarray*}
Then we get a contradiction with  \eqref{cc}, and since $k_0$ was arbitrary,
$\nu_{k}= 0$ for all $k\in I$. Hence as a consequence, $u_{n}\to u_0$ in
$L^{2^*_{\alpha}}(\Omega)$. We finish in the standard way: convergence of
$u_{n}$ in $L^{\frac{2N}{N-\alpha}}(\Omega)$ implies convergence of $f(u_{n})$
in $L^{\frac{2N}{N+\alpha}}(\Omega)$, and finally by using the continuity of the inverse operator
$(-\Delta)^{-\alpha/2}$, we obtain convergence of $u_{n}$ in
$H^{\alpha/2}_{0}(\Omega)$. ~\hfill$\Box$

\

Now it remains to show that we can obtain a local (PS)$_{c}$
sequence for $\widetilde{J}$ under the critical level $c=c^*$. To do that
  we will use $w_\varepsilon =\ext
(u_\varepsilon)$, the family of minimizers to the Trace inequality
\eqref{eq:trace-p}, where $u_\varepsilon$ is given in
\eqref{minimizers0}. We remark that, despite the cases $\alpha=1$
and $\alpha=2$,  $w_{\varepsilon}$ does not possesses an explicit expression.
This is an extra difficulty that we have to overcome. Taking into account that
the family $u_\varepsilon$ is self-similar,
$u_\varepsilon(x)=\varepsilon^{\frac{\alpha-N}{2}} u_1 ( x/\varepsilon)$ and
the fact that the Poisson kernel (\ref{poisson}) is also self-similar
\begin{equation}\label{nucleo_autosem}
P_{y}^{\alpha}(x)=\frac{1}{y^{N}}P^{\alpha}_1\left(\frac{x}{y}\right) ,
\end{equation}
gives easily that the family $w_\varepsilon$ satisfies
\begin{equation}\label{we}
w_{\varepsilon}(x,y)=\varepsilon^{\frac{\alpha-N}{2}}w_{1}\left(\frac{x}{\varepsilon},\frac{y}{\varepsilon}\right).
\end{equation}
We will denote $P^{\alpha}=P^{\alpha}_1$. Also, we will write
$w_{1,\alpha}$ instead of $w_1$ to emphasize the dependence on the
parameter $\alpha$.
\begin{Lemma}\label{propiedades_minimizantes}
With the above notation it holds
\begin{equation}\label{paso2}
|\nabla w_{1,\alpha}(x,y)|\leq \frac{C}{y}\,
w_{1,\alpha}(x,y),\quad\alpha>0,\: (x,y)\in\mathbb{R}^{N+1}_+
\end{equation}
and
\begin{equation}\label{paso3}
|\nabla w_{1,\alpha}(x,y)|\leq C
w_{1,\alpha-1}(x,y),\quad\alpha>1,\: (x,y)\in\mathbb{R}^{N+1}_+.
\end{equation}
\end{Lemma}
\pf Differentiating with  respect to each variable
$x_{i}\,,\,i=1,\,\ldots\,,N,$ and the variable $y$, it follows that
\begin{eqnarray*}
|\partial_{x_{i}}
w_{1,\alpha}(x,y)|&\leq&\int_{\mathbb{R}^{N}}{\frac{(N+\alpha)y^{\alpha}|x-z|}{(y^{2}+|x-z|^{2})^{\frac{N+\alpha}{2}+1}(1+|z|^{2})^{\frac{N-\alpha}{2}}}dz}\\
&\leq&\frac{N+\alpha}{2y}\int_{\mathbb{R}^{N}}{\frac{y^{\alpha}}{(y^{2}+|x-z|^{2})^{\frac{N+\alpha}{2}}(1+|z|^{2})^{\frac{N-\alpha}{2}}}dz}\\
&=&\frac{C}{y}w_{1,\alpha}(x,y)
\end{eqnarray*}
and
\begin{eqnarray*}
|\partial_{y}
w_{1,\alpha}(x,y)|&=& \left|{\int_{\mathbb{R}^{N}}{\frac{y^{\alpha-1}(\alpha|x-z|^{2}-Ny^{2})}{(y^{2}+|x-z|^{2})^{\frac{N+\alpha}{2}+1}(1+|z|^{2})^{\frac{N-\alpha}{2}}}dz}}\right|\\
&\leq&C\int_{\mathbb{R}^{N}}{\frac{y^{\alpha-1}}{(y^{2}+|x-z|^{2})^{\frac{N+\alpha}{2}}(1+|z|^{2})^{\frac{N-\alpha}{2}}}dz}\\
&=&\frac{C}{y}w_{1,\alpha}(x,y).
\end{eqnarray*}
Therefore we get (\ref{paso2}). To obtain (\ref{paso3}) we recall that
$u_{1,\alpha}(z)=(1+|z|^{2})^{-\frac{N-\alpha}{2}}.$ Then, by
(\ref{nucleo_autosem}) it follows that
\begin{eqnarray*}
\left|{\partial_{y}w_{1,\alpha}(x,y)}\right|&=&\left|{\partial_{y}\left(\int_{\mathbb{R}^{N}}{\frac{1}{y^{N}}P^{\alpha}\left(\frac{x-z}{y}\right)u_{1,\alpha}(z)dz}\right)}\right|\\
&=&\left|{-\partial_{y}\left(\int_{\mathbb{R}^{N}}{P^{\alpha}(\tilde{z})u_{1,\alpha}(x-y\tilde{z})d\tilde{z}}\right)}\right|\\
&=&\left|{\int_{\mathbb{R}^{N}}{P^{\alpha}(\tilde{z})\langle\tilde{z},\nabla u_{1,\alpha}(x-y\tilde{z})\rangle d\tilde{z}}}\right|\\
&=&\left|{-\int_{\mathbb{R}^{N}}{\frac{1}{y^{N}}P^{\alpha}\left(\frac{x-z}{y}\right)\langle\frac{x-z}{y},\nabla u_{1,\alpha}(z)\rangle dz}}\right|\\
&\leq&(N-\alpha)\int_{\mathbb{R}^{N}}{\frac{1}{y^{N}}P^{\alpha}\left(\frac{x-z}{y}\right)\frac{|x-z|}{y}\frac{|z|}{(1+|z|^{2})^{\frac{N-\alpha}{2}+1}}dz}\\
&\leq&(N-\alpha)\int_{\mathbb{R}^{N}}{\frac{y^{\alpha-1}}{(y^{2}+|x-z|^{2})^{\frac{N+\alpha-1}{2}}(1+|z|^{2})^{\frac{N-\alpha+1}{2}}}dz}\\
&=&Cw_{1,\alpha-1}(x,y).
\end{eqnarray*}
Doing the same calculations in variables $x_{i}$ for $i=1,\,\ldots\,,N$, we
obtain
\begin{eqnarray*}
\left|{\partial_{x_{i}}w_{1,\alpha}(x,y)}\right|&=&\left|{-\partial_{x_{i}}\left(\int_{\mathbb{R}^{N}}{P^{\alpha}(\tilde{z})u_{1,\alpha}(x-y\tilde{z})d\tilde{z}}\right)}\right|\\
&\leq&\int_{\mathbb{R}^{N}}{P^{\alpha}(\tilde{z})|\nabla u_{1,\alpha}|(x-y\tilde{z})d\tilde{z}}\\
&=&\int_{\mathbb{R}^{N}}{\frac{1}{y^{N}}P^{\alpha}\left(\frac{x-z}{y}\right)|\nabla u_{1,\alpha}|(z)dz}\\
&\leq&(N-\alpha)\int_{\mathbb{R}^{N}}{\frac{y^{\alpha}}{(y^{2}+|x-z|^{2})^{\frac{N+\alpha}{2}}}\frac{|z|}{(1+|z|^{2})^{\frac{N-\alpha}{2}+1}}dz}\\
&=&Cw_{1,\alpha-1}(x,y).
\end{eqnarray*}
~\hfill$\Box$
\\
Let us now introduce a cut-off function $\phi_{0}(s) \in
C^{\infty}(\mathbb{R}_{+})$, nonincreasing satisfying
$$
\phi_{0}(s) = 1 \mbox{ if } 0\leq s \leq \frac{1}{2},\quad
\phi_{0}(s) = 0 \mbox{ if } s\geq 1.
$$
Assume without loss of generality that $0\in \Omega$. We then define, for some
fixed $r>0$ small enough such that
$\overline{B}^{+}_{r}\subseteq\overline{\mathcal{C}}_{\Omega}$, the function
$\phi(x,y)=\phi_{r}(x,y)=\phi_{0}(\frac{r_{xy}}{r})$ with
$r_{xy}=|(x,y)|=(|x|^2+y^2)^{1/2}$.  Note that $\phi\omega_{\varepsilon}\in
X^{\alpha}_{0}(\mathcal{C}_{\Omega})$. Thus we get
\begin{Lemma}\label{propiedades_eta}
With the above notation, the family $\{\phi w_\varepsilon\}$, and its trace on
$\{y=0\}$, namely $\{\phi u_\varepsilon\}$, satisfy
\begin{equation}\label{norm-estimates}
\|\phi
w_{\varepsilon}\|_{X^{\alpha}_{0}(\mathcal{C}_{\Omega})}^{2}=
\|w_{\varepsilon}\|_{X^{\alpha}_{0}(\mathcal{C}_{\Omega})}^{2} +
O(\varepsilon^{N-\alpha}),
\end{equation}
\begin{equation}\label{norm-estimates-sig}
\|\phi u_{\varepsilon}\|_{L^{2}(\Omega)}^{2}=\left\{
\begin{array}{ll}
C\varepsilon^{\alpha}+ O(\varepsilon^{N-\alpha}) &\mbox{ if } N> 2\alpha,\\
C\varepsilon^{\alpha}{\color{black}\log(1\slash \varepsilon)}+
O(\varepsilon^{\alpha}) &\mbox{ if }  N=2\alpha,
\end{array}\right.
\end{equation}
and
\begin{equation}\label{norm-r-estimates}
\|\phi
u_{\varepsilon}\|_{L^r(\Omega)}^{r}\ge c\varepsilon^{\frac{N-\alpha}{2}},\quad \alpha<N<2\alpha,\quad r=\frac{N+\alpha}{N-\alpha},
\end{equation}
 for $\varepsilon$ small enough and $C>0$.
\end{Lemma}
\pf
The product $\phi w_{\varepsilon}$ satisfies
\begin{eqnarray}\label{norm-estimates2}
\|\phi w_{\varepsilon}\|_{X^{\alpha}_{0}(\mathcal{C}_{\Omega})}^{2}&=&\kappa_{\alpha}\int_{\mathcal{C}_{\Omega}}{y^{1-\alpha}(|\phi
\nabla w_{\varepsilon}|^2+|w_{\varepsilon}\nabla
\phi |^2+2\langle w_{\varepsilon} \nabla\phi ,\phi \nabla
w_{\varepsilon}
\rangle)dxdy}\nonumber\\
&\leq&\|w_{\varepsilon}\|_{X^{\alpha}_{0}(\mathcal{C}_{\Omega})}^{2}
+\kappa_{\alpha}\int_{\mathcal{C}_{\Omega}}y^{1-\alpha}|w_{\varepsilon}
\nabla \phi |^2dxdy  \\
& & + 2\kappa_{\alpha}
\int_{\mathcal{C}_{\Omega}}y^{1-\alpha}\langle w_{\varepsilon}
\nabla\phi ,\phi  \nabla
w_{\varepsilon}\rangle dxdy.\nonumber
\end{eqnarray}
To estimate the second term of the right hand side, we observe that $0\leq u_{\varepsilon}(x) \leq
\varepsilon^{\frac{N-\alpha}{2}} |x|^{\alpha-N} $, and since the extension of
the function $\Gamma(x)=|x|^{\alpha-N}$ is $\widetilde\Gamma(x,y)=(|x|^2 +
y^2)^{\frac{\alpha-N}2}= r^{\alpha-N}_{xy}$, we get that
\begin{eqnarray}
\int_{\mathcal{C}_{\Omega}}{y^{1-\alpha}|w_{\varepsilon} \nabla
\phi |^2dxdy}&\leq&C\int_{\{\frac{r}{2}\leq
r_{xy} \leq r\}}{y^{1-\alpha}w_{\varepsilon}^2dxdy}\nonumber\\
&\leq&C\varepsilon^{N-\alpha}\int_{\{\frac{r}{2}\leq r_{xy} \leq
r\}}{y^{1-\alpha}
r_{xy}^{2(\alpha-N)}dxdy}\\
&=&O(\varepsilon^{N-\alpha})\nonumber.\label{termino2}
\end{eqnarray}
For the remaining term we need to use the properties of the function
$w_{\varepsilon}$ given in Proposition
\ref{propiedades_minimizantes}. Let $C_{r}=\{r/2\leq r_{xy}\leq r\}
\subset \mathcal{C}_{\Omega}$. By (\ref{we}) we get
\begin{eqnarray}\label{homog}
\int_{\mathcal{C}_{\Omega}}{y^{1-\alpha}\langle
w_{\varepsilon}\nabla\phi ,\phi \nabla w_{\varepsilon}\rangle
dxdy}&\leq&C\int_{C_{r}}{y^{1-\alpha}|w_{\varepsilon}(x,y)\|\nabla
w_{\varepsilon}(x,y)|dxdy}\nonumber\\
&=&C\varepsilon^{-N+\alpha-1}\int_{C_{r}}{y^{1-\alpha}\left|{w_{1,\alpha}\left(\frac{x}{\varepsilon},\frac{y}{\varepsilon}\right)}
\right|\left|{\nabla
w_{1,\alpha}\left(\frac{x}{\varepsilon},\frac{y}{\varepsilon}\right)}\right|dxdy}\\
&=&
C\varepsilon\int_{C_{\frac{r}{\varepsilon}}}{y^{1-\alpha}|w_{1,\alpha}(x,y)|\,|\nabla
w_{1,\alpha}(x,y)|dxdy}\nonumber.
\end{eqnarray}
Moreover, for $(x,y)\in C_{r\slash\varepsilon}$ and $\alpha>0$, we
obtain that
\begin{eqnarray}\label{blabla}
w_{1,\alpha}(x,y)&=&\int_{|z|<\frac{1}{4\varepsilon}}{P_{y}^{\alpha}(x-z)u_{1,\alpha}(z)dz}+\int_{|z|>\frac{1}{4\varepsilon}}{P_{y}^{\alpha}(x-z)u_{1,\alpha}(z)dz}\nonumber\\
&\leq&C\varepsilon^{N+\alpha}y^{\alpha}\int_{|z|<\frac{1}{4\varepsilon}}{\frac{dz}{|z|^{N-\alpha}}}+C\varepsilon^{N-\alpha}\int_{\mathbb{R}^{N}}{P_{y}^{\alpha}(z)dz}\\
&\leq&Cy^{\alpha}\varepsilon^{N}+C\varepsilon^{N-\alpha}\nonumber \leq C\varepsilon^{N-\alpha}.
\end{eqnarray}
If $\alpha<1$, from (\ref{paso2}), (\ref{homog}) and (\ref{blabla}), it follows
that
\begin{equation}\label{menor_que_1}
\int_{\mathcal{C}_{\Omega}}{y^{1-\alpha}\langle
w_{\varepsilon}\nabla\phi ,\phi \nabla w_{\varepsilon}\rangle dxdy}\leq
C\varepsilon^{1+2(N-\alpha)}\int_{C_{\frac{r}{\varepsilon}}}{y^{-\alpha}dxdy}=O(\varepsilon^{N-\alpha}).
\end{equation}
To obtain the similar estimate for $\alpha>1$ we use (\ref{paso3}).
Indeed by this estimate, together with (\ref{homog}) and
(\ref{blabla}) we get that
\begin{equation}\label{mayor_que_1}
\int_{\mathcal{C}_{\Omega}}{y^{1-\alpha}\langle
w_{\varepsilon}\nabla\phi ,\phi \nabla w_{\varepsilon}\rangle dxdy}\leq
C\varepsilon^{2(1+N-\alpha)}\int_{C_{\frac{r}{\varepsilon}}}{y^{1-\alpha}dxdy}=O(\varepsilon^{N-\alpha}).
\end{equation}
Note that for $\alpha = 1$, as $w_{\varepsilon}$ is explicit, we can
obtain the same estimate directly.

Then we have proved that
\[
\|\phi
w_{\varepsilon}\|_{X^{\alpha}_{0}(\mathcal{C}_{\Omega})}^{2}=
\|w_{\varepsilon}\|_{X^{\alpha}_{0}(\mathcal{C}_{\Omega})}^{2} +
O(\varepsilon^{N-\alpha}).
\]

We now show that (\ref{norm-estimates-sig}) holds. 
\begin{eqnarray*}
\|\phi
u_{\varepsilon}\|_{L^{2}(\Omega)}^{2}&=&\int_{\Omega}{\phi ^{2}(x)
\frac{\varepsilon^{N-\alpha}}{(|x|^2+\varepsilon^2)^{N-\alpha}}dx}\\
&\geq&\int_{\{|x|<r\slash2\}}{\frac{\varepsilon^{N-\alpha}}{(|x|^2+\varepsilon^2)^{N-\alpha}}dx}\\
&\geq&\int_{\{|x|<\varepsilon\}}{\frac{\varepsilon^{N-\alpha}}{(2\varepsilon^2)^{N-\alpha}}dx}+
\int_{\{\varepsilon<|x|<r\slash2\}}{\frac{\varepsilon^{N-\alpha}}{(2|x|^2)^{N-\alpha}}dx}\\
&=&C\varepsilon^{\alpha}+C\varepsilon^{N-\alpha}\int^{r/2}_{\varepsilon}{\theta^{2\alpha
- 1 -N}d\theta.}
\end{eqnarray*}
Finally, \eqref{norm-r-estimates} follows in a similar way to \eqref{norm-estimates-sig}, so we omit the details.
~\hfill$\Box$

With the above properties in mind, we define the family of functions
$\eta_{\varepsilon}= \frac{\phi  w_{\varepsilon}}{\|\phi
u_{\varepsilon}\|_{L^{2^*_{\alpha}}(\Omega)}}$.
\begin{Lemma}\label{underlevel}
There exists $\varepsilon
>0$ small enough such that
\begin{equation}\label{lowlevel}
\displaystyle \sup_{t\geq 0}\widetilde{J}(t\eta_{\varepsilon}) <
\textcolor{black}{c^*}.
\end{equation}
\end{Lemma}
\pf Assume $N\ge 2\alpha$, we make use of the following estimate
\begin{equation}\label{estimate-ABC1}
(a+b)^p\geq a^p + b^p + \mu a^{p-1}b, \quad a,b \geq 0, \quad p>1,\quad\mbox{for some}\, \mu>0.
\end{equation}
Therefore
\begin{equation}\label{estimatesG}
G(w) \geq \frac{1}{2^*_{\alpha}}w^{2^*_{\alpha}}+\frac{\mu}{2}w^2
w_{0}^{2^*_{\alpha} -2}
\end{equation}
which implies
$$
\widetilde{J}(t\eta_{\varepsilon}) \leq
\frac{t^2}{2}\|\eta_{\varepsilon}\|_{X^{\alpha}_{0}(\mathcal{C}_{\Omega})}^{2}
- \frac{t^{2^*_{\alpha}}}{2^*_{\alpha}}-\frac{t^2}{2}\mu
\int_{\Omega}w_{0}^{2^*_{\alpha}-2}\eta_{\varepsilon}^2dx.
$$
Since $w_{0}\ge a_0>0$ in $\supp(\eta_{\varepsilon})$ we get
$$
\widetilde{J}(t\eta_{\varepsilon}) \leq
\frac{t^2}{2}\|\eta_{\varepsilon}\|_{X^{\alpha}_{0}(\mathcal{C}_{\Omega})}^{2}
-
\frac{t^{2^*_{\alpha}}}{2^*_{\alpha}}-\frac{t^2}{2}\widetilde{\mu}
\|\eta_{\varepsilon}\|_{L^{2}(\Omega)}^{2} =:g(t).
$$
It is clear that $\displaystyle \lim_{t\to\infty}
g(t) = -\infty$, and
$\displaystyle\sup_{t\geq 0}
g(t)$ is attained at some
$t_{\varepsilon} > 0$. By differentiating the above function we obtain
\begin{equation}\label{difJ}
0=g'(t_{\varepsilon})=
t_{\varepsilon}\|\eta_{\varepsilon}\|_{X^{\alpha}_{0}(\mathcal{C}_{\Omega})}^{2}-t_{\varepsilon}^{2^*_{\alpha}
-1}-t_{\varepsilon}\widetilde{\mu}\|\eta_{\varepsilon}\|_{L^{2}(\Omega)}^{2},
\end{equation}
which implies
$$
t_{\varepsilon} \leq
\|\eta_{\varepsilon}\|_{X^{\alpha}_{0}(\mathcal{C}_{\Omega})}^{\frac{2}{2^*_{\alpha}-2}}.
$$
Observe that by Lemma \ref{propiedades_eta} we have $t_\varepsilon\ge C>0$. On
the other hand, the function
$$
t \mapsto
\frac{t^2}{2}\|\eta_{\varepsilon}\|_{X^{\alpha}_{0}(\mathcal{C}_{\Omega})}^{2}
- \frac{t^{2^*_{\alpha}}}{2^*_{\alpha}}
$$
is increasing on
$[0,\|\eta_{\varepsilon}\|_{X^{\alpha}_{0}(\mathcal{C}_{\Omega})}^{\frac{2}{2^*_{\alpha}-2}}]$.
Whence
$$
\displaystyle \sup_{t\geq 0}
g(t)=g(t_{\varepsilon})\leq
\frac{\alpha}{2N}\|\eta_{\varepsilon}\|_{X^{\alpha}_{0}(\mathcal{C}_{\Omega})}^{\frac{2N}{\alpha}}-
C\|\eta_{\varepsilon}\|_{L^{2}(\Omega)}^{2}.
$$
Since $\|u_{\varepsilon}\|_{L^{2^*_{\alpha}}(\Omega)}$ is independent of $\varepsilon$,
by Lemma \ref{propiedades_eta}  we have
\begin{equation}\label{estimatesOld}
\|\eta_{\varepsilon}\|_{X^{\alpha}_{0}(\mathcal{C}_{\Omega})}^{2}=
\textcolor{black}{\kappa_{\alpha}}S(\alpha,N) +
O(\varepsilon^{N-\alpha})
\end{equation}
and
$$
\|\eta_{\varepsilon}\|_{L^{2}(\Omega)}^{2} = \left\{
\begin{array}{ll}
O(\varepsilon^{\alpha}) &\mbox{ if } N> 2\alpha,\\
O(\varepsilon^{\alpha} \log(1\slash \varepsilon)) &\mbox{ if }
N=2\alpha.
\end{array}\right.
$$
Therefore, for $N>2\alpha$, we get that
\begin{equation}\label{estimateJ}
g(t_{\varepsilon})\leq\frac{\alpha}{2N}(\textcolor{black}{\kappa_{\alpha}}S(\alpha,N))^{\frac{N}{\alpha}}+C\varepsilon^{N-\alpha}-C\varepsilon^{\alpha}<\frac{\alpha}{2N}(\textcolor{black}{\kappa_{\alpha}}S(\alpha,N))^{\frac{N}{\alpha}}=\textcolor{black}{c^*}.
\end{equation}
If $N=2\alpha$ the same conclusion follows.

The last case $\alpha <N<2\alpha$ follows by using the estimate \eqref{estimate-ABC1} which gives
\begin{equation}\label{estimatesG2}
G(w) \geq \frac{1}{2^*_{\alpha}}w^{2^*_{\alpha}}+w_0 w^{2^*_{\alpha} -1}.
\end{equation}
Then \eqref{estimatesG2} jointly with \eqref{norm-r-estimates} and arguing in a similar way as above finish the proof.
~\hfill$\Box$

\pff{Theorem~{\rm \ref{maintheorem2}}-\rm(3)}

To finish the last statement
in Theorem \ref{maintheorem2}, in view of the previous results, we seek for critical values
below level $c^*$. For that purpose, we want to use  the classical MP Theorem by Ambrosetti-Rabinowitz in \cite{AR}.
We define
$$
\Gamma_\varepsilon=\{ \gamma\in\mathcal{C}([0,1],X_0^\alpha(\mathcal{C}_\Omega)) :\: \gamma(0)=0,\, \gamma(1)=t_\varepsilon \eta_\varepsilon\}
$$
for some $t_\varepsilon>0 $ such that $\widetilde{J}(t_\varepsilon \eta_\varepsilon)<0$. And consider the minimax value
$$
c_\varepsilon=\inf_{\gamma\in \Gamma_\varepsilon}\max \{\widetilde{J}(\gamma(t)):\: 0\le t\le 1\}.$$
According to Lemma  \ref{lemma:minimo}, $c_\varepsilon\ge 0$. By Lemma \ref{underlevel},  for $\varepsilon\ll 1$,
$$c_\varepsilon\le \sup_{t\ge 0}\widetilde{J}(t \eta_\varepsilon)<c^*=\frac{\alpha}{2N}\left(\kappa_\alpha S(\alpha,N)\right)^{N/\alpha}.$$
This estimate jointly with Lemma \ref{strongly_convergent} and the MPT \cite{AR} if the minimax energy level is positive, or the refinement of the MPT \cite{GhP} if the minimax level is zero, give the existence of a second solution to $(P)_\lambda$. ~\hfill$\Box$


\section{Linear and superlinear cases.}\label{sect-linear-superlinear}
\subsection{Linear case}
\setcounter{equation}{0} The proof of Theorem \ref{maintheorem}
follows the ideas of \cite{brezis-nirenberg}. Note that for
$\alpha=1$, where the minimizers given in (\ref{we}) are explicit,
this result was recently proved in \cite{tan}.
\\The first part of that theorem is an
straightforward calculus.

\pff{Theorem~{\rm \ref{maintheorem}} \rm(1)}
Let $\varphi_{1}$ be the first eigenfunction of
$(-\Delta)^{\alpha/2}$ in $\Omega$. We have
$$
\int_{\Omega} (-\Delta)^{\alpha/4}u(-\Delta)^{\alpha/4}\varphi_1\,
dx = \int_{\Omega} \lambda_{1} u\varphi_{1} \, dx.
$$
On the other hand,
$$
\int_{\Omega} (-\Delta)^{\alpha/4}u(-\Delta)^{\alpha/4}\varphi_1\,
dx = \int_{\Omega}[u^{2^*_{\alpha}-1} + \lambda u] \varphi_{1}\,
dx > \int_{\Omega} \lambda u \varphi_{1}\, dx .
$$
This clearly implies
$\lambda
<\lambda_1$.
~\hfill$\Box$

To prove the second part of Theorem \ref{maintheorem}  some notation is in order. We consider the following Rayleigh
quotient
$$
    Q_{\lambda}(w)= \frac{\|w\|_{X^{\alpha}_{0}(\mathcal{C}_{\Omega})}^{2} - \lambda \|u\|_{L^{2}
    (\Omega)}^{2}}{ \|u\|_{L^{2^*_{\alpha}}(\Omega)}^{2}}
$$
and
    \begin{equation}\label{Slambda}
    S_{\lambda}= \inff\{Q_{\lambda}(w)\, |\: \ w\in X^{\alpha}_{0}(\mathcal{C}_{\Omega})
    \}.
    \end{equation}

\begin{Proposition}\label{inequality}
Assume $0< \lambda <\lambda_{1}$. Then
$S_{\lambda}<\kappa_{\alpha}S(\alpha , N).$
\end{Proposition}
\pf Let $\phi=\phi_r $ be a cut-off function like in Lemma
\ref{propiedades_eta} and denote $\phi(x):=\phi (x,0)$. Taking
$r$ sufficiently small we can use $\phi w_{\varepsilon}\in
X^{\alpha}_{0}(\mathcal{C}_{\Omega})$ as a test function in $Q_{\lambda}$,
where $w_{\varepsilon}$ is defined in (\ref{we}).
Denoting $K_1=\|u_{\varepsilon}\|_{L^{2^*_{\alpha}}(\Omega)}^{2^*_\alpha}$, as before, $K_1$ is independent of $\varepsilon$, and moreover
\begin{eqnarray}\label{primera_parte}
\int_{\Omega}{|\phi u_{\varepsilon}|^{2^*_{\alpha}}dx}&=&\int_{\mathbb{R}^{N}}{|\phi u_{\varepsilon}|^{2^*_{\alpha}}dx}\nonumber\\
&\geq&\int_{|x|<r/2}{|u_{\varepsilon}|^{2^*_{\alpha}}dx}\nonumber\\
&=&K_{1}-\int_{|x|>r/2}{|u_{\varepsilon}|^{2^*_{\alpha}}dx}\nonumber\\
&\geq&K_{1} + O(\varepsilon^{N}).
\end{eqnarray}
Since $w_{\varepsilon}$ is a minimizer of $S(\alpha , N)$, we have
that
\begin{equation}\label{segunda_parte}
K_{1}^{-2\slash 2^*_{\alpha}} \int_{\mathbb{R}^{N+1}_{+}}
y^{1-\alpha}|\nabla w_{\varepsilon}|^{2} \, dxdy =S(\alpha ,
N).
\end{equation}
\\Finally, by (\ref{primera_parte}) and using the estimates (\ref{norm-estimates}) and (\ref{norm-estimates-sig}), for $N>2\alpha$, we obtain that
$$
Q_{\lambda}(\phi  w_{\varepsilon})\leq \frac{\displaystyle
\kappa_{\alpha}\int_{\mathbb{R}^{N+1}_{+}} y^{1-\alpha}|\nabla
w_{\varepsilon}|^{2} \, dxdy-\lambda C\varepsilon^{\alpha} +
O({\varepsilon^{N-\alpha}})}{K_{1}^{2\slash 2^*_{\alpha}} +
O(\varepsilon^{N})}.
$$
Therefore taking $\varepsilon$ small enough, we get
\begin{eqnarray*}
\displaystyle Q_{\lambda}(\phi  w_{\varepsilon})&\leq&
\frac{\kappa_{\alpha}S(\alpha,N) - \lambda C\varepsilon^{\alpha}
K_{1}^{-2\slash 2^*_{\alpha}} +
O(\varepsilon^{N-\alpha})}{1+O(\varepsilon^{N})} \\[3mm]
\displaystyle &\leq&  \kappa_{\alpha}S(\alpha,N) - \lambda
C\varepsilon^{\alpha} K_{1}^{-2\slash 2^*_{\alpha}} +
O(\varepsilon^{N-\alpha}) \\[3mm]
\displaystyle &<& \kappa_{\alpha}S(\alpha,N).
\end{eqnarray*}

On the other hand, a similar calculus for the case $N=2\alpha$,
proves that for $\varepsilon$ small enough,
$$
Q_{\lambda}(\phi  w_{\varepsilon})\leq
\kappa_{\alpha}S(\alpha,N) - \lambda C\varepsilon^{\alpha}
\log(1\slash \varepsilon) K_{1}^{-2\slash 2^*_{\alpha}} +
O(\varepsilon^{\alpha})< \kappa_{\alpha}S(\alpha,N),
$$
which finishes the proof. ~\hfill$\Box$

Recall now the Brezis-Lieb Lemma,
\begin{Lemma}[\cite{bl}]\label{brezis_lieb}
Let $\Omega$ be an open set and $\{u_{n}\}$ be a sequence weakly convergent in $L^{q}(\Omega), 2\leq
q<\infty$ and a.e. convergent in $\Omega$. Then
$\displaystyle\lim_{n\to\infty}(\|u_{n}\|^{q}_{L^{q}(\Omega)}-\|u_{n}-u\|^{q}_{L^{q}(\Omega)})=\|u\|^{q}_{L^{q}(\Omega)}.$
\end{Lemma}
This property allows us to we prove the following one.
\begin{Proposition}\label{achieved}
Assume $0<\lambda<\lambda_{1}$. Then the infimum $S_{\lambda}$ defined  in
\eqref{Slambda} is achieved.
\end{Proposition}
\pf First, since $\lambda<\lambda_{1}$ we have that $S_{\lambda}>0$.
Let us take a minimizing sequence of $S_{\lambda}$, $\{w_{m}\}\subset
X^{\alpha}_{0}(\mathcal{C}_{\Omega}) $  such that,
without loss of generality, $w_{m}\geq 0$ and $\|w_{m}(\cdot,
0)\|_{L^{2^*_{\alpha}}(\Omega)}=1$. Clearly this implies that
$\|w_{m}\|_{X^{\alpha}_{0}(\mathcal{C}_{\Omega})} \leq C$, then there exists a subsequence
(still denoted by $\{w_{m}\}$) verifying
\begin{equation} \label{conditions}
\begin{array}{rcll}
 w_{m} &\rightharpoonup& w \quad\qquad \mbox{weakly in } X^{\alpha}_{0}(\mathcal{C}_{\Omega}),\\
 w_{m}(\cdot,0) &\to& w(\cdot,0) \quad \mbox{strongly in } L^{q}(\Omega), \; 1\leq q < 2^*_{\alpha} ,\\
 w_{m}(\cdot,0) &\to& w(\cdot,0) \quad \mbox{a.e in }
 \Omega.
\end{array}
\end{equation}
A simple calculation, using the weak convergence, gives that
\begin{eqnarray*}
\|w_{m}\|_{X^{\alpha}_{0}(\mathcal{C}_{\Omega})}^{2}&=&
\|w_{m}-w\|_{X^{\alpha}_{0}(\mathcal{C}_{\Omega})}^{2} +
\|w\|_{X^{\alpha}_{0}(\mathcal{C}_{\Omega})}^{2} +
2\kappa_{\alpha}\int_{\mathcal{C}_{\Omega}}{y^{1-\alpha}\langle \nabla w, \nabla
w_{m}- \nabla w \rangle dxdy}\\
&=&\|w_{m}-w\|_{X^{\alpha}_{0}(\mathcal{C}_{\Omega})}^{2} +
\|w\|_{X^{\alpha}_{0}(\mathcal{C}_{\Omega})}^{2} + o(1).
\end{eqnarray*}
By Lemma \ref{brezis_lieb}, we have that $\|(w_{m}-w)(\cdot,
0)\|_{L^{2^*_{\alpha}}(\Omega)}\leq1$ for $m$ big enough. Hence
\begin{eqnarray*}
Q_{\lambda}(w_{m})&=&\|w_{m}\|_{X^{\alpha}_{0}(\mathcal{C}_{\Omega})}^{2}
- \lambda
\|w_{m}(\cdot, 0)\|_{L^{2}(\Omega)}^{2}\\
&=&\|w_{m}-w\|_{X^{\alpha}_{0}(\mathcal{C}_{\Omega})}^{2}+
\|w\|_{X^{\alpha}_{0}(\mathcal{C}_{\Omega})}^{2} - \lambda
\|w_{m}(\cdot, 0)\|_{L^{2}(\Omega)}^{2} + o(1)\\
&\geq& \kappa_{\alpha}S(\alpha,N)\|(w_{m}-w)(\cdot,
0)\|_{L^{2^*_{\alpha}}(\Omega)}^{2}+S_{\lambda}\|w(\cdot,
0)\|_{L^{2^*_{\alpha}}(\Omega)}^{2}
+ o(1)\\
&\geq& \kappa_{\alpha}S(\alpha,N)\|(w_{m}-w)(\cdot,
0)\|_{L^{2^*_{\alpha}}(\Omega)}^{2^*_{\alpha}}+S_{\lambda}\|w(\cdot,
0)\|_{L^{2^*_{\alpha}}(\Omega)}^{2^*_{\alpha}} + o(1).
\end{eqnarray*}
By Lemma \ref{brezis_lieb} again, this leads to
\begin{eqnarray*}
Q_{\lambda}(w_{m})&\geq&
(\kappa_{\alpha}S(\alpha,N)-S_{\lambda})\|(w_{m}-w)(\cdot,
0)\|_{L^{2^*_{\alpha}}(\Omega)}^{2^*_{\alpha}}+S_{\lambda}\|w_{m}(\cdot,
0)\|_{L^{2^*_{\alpha}}(\Omega)}^{2^*_{\alpha}} + o(1)\\
&=&(\kappa_{\alpha}S(\alpha,N)-S_{\lambda})\|(w_{m}-w)(\cdot,
0)\|_{L^{2^*_{\alpha}}(\Omega)}^{2^*_{\alpha}}+S_{\lambda} + o(1).
\end{eqnarray*}
Since $\{w_{m}\}$ is a minimizing sequence for $S_{\lambda}$, we
obtain:
$$
o(1)+S_{\lambda} \geq
(\kappa_{\alpha}S(\alpha,N)-S_{\lambda})\|(w_{m}-w)(\cdot,
0)\|_{L^{2^*_{\alpha}}(\Omega)}^{2^*_{\alpha}}+S_{\lambda} + o(1).
$$
Thus by Proposition \ref{inequality}
$$
w_{m}(\cdot,0) \to w(\cdot,0) \qquad \mbox{ in }
L^{2^*_{\alpha}}(\Omega).
$$
Finally, by a standard lower semi-continuity argument,  $w$
is a minimizer for $Q_{\lambda}$. ~\hfill$\Box$

\pff{Theorem~{\rm \ref{maintheorem}} \rm(2)}
By Proposition \ref{achieved} there exists an $\alpha$-harmonic function $w \in
X^{\alpha}_{0}(\mathcal{C}_{\Omega})$,  such that $
\|u\|_{L^{2^*_{\alpha}}(\Omega)}^{2}=1$ and
$$
\|w\|_{X^{\alpha}_{0}(\mathcal{C}_{\Omega})}^{2} - \lambda \|u\|_{L^{2}
(\Omega)}^{2} = S_{\lambda}
$$where $u=w(\cdot,0)$. Without loss of generality we may assume $w\geq 0$ (otherwise we take $|w|$ instead of $w$). So we get a positive solution of $(P_{\lambda})$.
~\hfill$\Box$
\subsection{Superlinear case.}

In order to prove Theorem \ref{maintheorem3}, the only difficult part is to
show that we have a (PS)$_c$ sequence under the critical level $c=c^*$. This
follows the same type of computations like in Lemma~\ref{underlevel}, with the estimate
$\|\eta_{\varepsilon}\|_{L^{q+1}(\Omega)}^{q+1}\ge C\varepsilon^{\frac{\alpha-N}{2}q+\frac{\alpha+N}{2}}$ which holds for $N>\alpha(1+\frac{1}{q})$. In this case there is no limitation on $\lambda>0$. We omit the complete details.

\section{Regularity \& Concentration-Compactness}\label{sect-regularity}
\setcounter{equation}{0}

We begin this section with some results about the boundedness and regularity of
solutions. The next proposition is a refinement of Proposition 5.3 of
\cite{brandle-colorado-depablo-sanchez} in order to cover the critical case
$p=2^{*}_{\alpha}-1$. It is essentially based on
\cite{brezis-kato}.
\begin{Proposition}\label{prop:bound}
Let $u \in H_{0}^{\alpha\slash 2}(\Omega)$ be a solution to the
problem
\begin{equation}\label{p}
   \left\{
 \begin{array}{ll}
    (-\Delta)^{\alpha/2}u=f(x,u)&\quad\mbox{\rm in } \Omega, \\
    u>0&\quad\mbox{\rm in } \Omega, \\
    u=0&\quad\mbox{\rm on } \partial\Omega
  \end{array}\right.
\end{equation}
with $f$ satisfying
\begin{equation}\label{growth_condition}
0\leq f(x,s) \leq C(1+|s|^p) \quad \forall\,
(x,s)\in\Omega\times\mathbb{R}, \mbox{ and some }  0<p\leq 2^{*}_{\alpha}-1.
\end{equation}
Then $u\in L^{\infty}(\Omega)$  with
$\|u\|_{L^{\infty}(\Omega)}\leq C(\|u\|_{H_{0}^{\alpha\slash
2}(\Omega)})$.
\end{Proposition}
\pf Let $w \in X^{\alpha}_{0}(\mathcal{C}_{\Omega})$ be a solution
to the problem
\begin{equation}\label{pp}
 \left\{
 \begin{array}{ll}
    L_{\alpha}w=0 &\quad \mbox{in } \mathcal{C}_{\Omega} , \\ [1mm]
    \dfrac{\partial w}{\partial \nu^{\alpha}}=f(\cdot,w)
    &\quad \mbox{in } \Omega, \\ [2mm]
    w=0 & \quad \mbox{on } \partial_{L} \mathcal{C}_{\Omega}.
  \end{array}
  \right.
\end{equation}
Then $u=w(\cdot,0)$ is a solution to (\ref{p}). Let
$$
a(x):=\frac{f(x,u)}{1+u(x)}.
$$
Clearly
\begin{equation}\label{a}
0\leq a \leq C(1+u^{p-1}) \in L^{\frac{N}{\alpha}}(\Omega),\quad\mbox{for}\,0<p\leq 2^*_{\alpha}-1.
\end{equation}
Given $T>0$ we denote
$$
w_{T}=w-(w-T)_{+}, \quad u_{T}=w_{T}(\cdot,0).
$$
For $\beta \geq 0$ we have
$$
\begin{array}{rl}
\displaystyle
\|ww_{T}^{\beta}\|_{X^{\alpha}_0(C_{\Omega})}^{2}=
&\displaystyle\kappa_{\alpha}\int_{\mathcal{C}_{\Omega}}y^{1-\alpha}w_{T}^{2\beta}|\nabla
w|^2 \,dxdy   \\
&\displaystyle  + \kappa_{\alpha}(2\beta + \beta^2)\int_{\{w\leq \,
T\}}y^{1-\alpha}w^{2\beta}|\nabla w|^2 \,dxdy.
\end{array}
$$
Using $\varphi=ww^{2\beta}_{T}\in X^{\alpha}_0(C_{\Omega})$ as a test function
we obtain
$$
\kappa_{\alpha}\int_{\mathcal{C}_{\Omega}} y^{1-\alpha} \langle\nabla w,
 \nabla (ww^{2\beta}_{T})\rangle\, dxdy = \int_{\Omega}
f(u)uu^{2\beta}_{T}\, dx \leq 2\int_{\Omega}
a(1+u^2)u^{2\beta}_{T}\, dx.
$$
On the other hand, it is clear that
$$
\begin{array}{ll}
\displaystyle \int_{\mathcal{C}_{\Omega}} y^{1-\alpha} \langle\nabla w,
 \nabla (ww^{2\beta}_{T})\rangle\, dxdy
=\int_{\mathcal{C}_{\Omega}}y^{1-\alpha}w_{T}^{2\beta}|\nabla w|^2 \,dxdy+ \\
[3mm]
\displaystyle + 2\beta\int_{\{w\leq \,
T\}}y^{1-\alpha}w^{2\beta}|\nabla w|^2 \,dxdy.
\end{array}
$$
Summing up, we have
$$
\|ww_{T}^{\beta}\|_{X^{\alpha}_0(C_{\Omega})}^{2}\leq C
\int_{\Omega} a(1+u^2)u^{2\beta}_{T}\, dx,
$$
which by  \eqref{eq:trace-r} implies that
\begin{equation}\label{traz}
\|uu_{T}^{\beta}\|_{L^{2^*_{\alpha}}(\Omega)}^{2}\leq \widetilde{C}
\int_{\Omega} a(1+u^2)u^{2\beta}_{T}\, dx,
\end{equation}
with $\widetilde{C}$ some positive constant depending on $\alpha$, $\beta$,
$N$ and $|\Omega|$. To compute the term on the right-hand side we
add the hypothesis $u^{\beta+1}\in L^{2}(\Omega)$. With this
assumption we get
$$
\begin{array}{lc}
 \displaystyle \int_{\Omega} a u^2u^{2\beta}_{T}\, dx \leq T_{0}\int_{\{a< \,
T_{0}\}} u^2u^{2\beta}_{T}\, dx + \int_{\{a\geq \, T_{0}\}}
au^2u^{2\beta}_{T}\, dx \\[3mm]
\displaystyle \leq C_{1}T_{0} + \left(\int_{\{a\geq \, T_{0}\}}
a^{\frac{N}{\alpha}} \,dx
\right)^{\frac{\alpha}{N}}\left(\int_{\Omega}
(uu^{\beta}_{T})^{2^*_{\alpha}}\, dx
\right)^{\frac{2}{2^*_{\alpha}}}.
\end{array}
$$
By the same calculation,
$$
\displaystyle \int_{\Omega} a u^{2\beta}_{T}\, dx \leq
C_{2}T_{0} + \left(\int_{\{a\geq \, T_{0}\}} a^{\frac{N}{\alpha}}
\,dx \right)^{\frac{\alpha}{N}}\left(\int_{\Omega}
(u^{\beta}_{T})^{2^*_{\alpha}}\, dx
\right)^{\frac{2}{2^*_{\alpha}}},
$$
where, since $u^{\beta+1}\in L^{2}(\Omega)$, $C_{1}$ and $C_{2}$
can be taken independent of $T$. Hence, by (\ref{a}) it follows
that
$$
\epsilon(T_{0})=\left(\int_{\{a\geq \, T_{0}\}}
a^{\frac{N}{\alpha}} \,dx \right)^{\frac{\alpha}{N}}\to 0
\quad \mbox{ as }\, T_{0} \to \infty.
$$
Therefore, choosing $T_{0}$ large enough such that
$C\epsilon(T_{0})<\frac{1}{2}$, by (\ref{traz}), we obtain that
there exists a constant $K(T_{0})$, independent of $T$, for which it
holds
$$
\|uu_{T}^{\beta}\|_{L^{2^*_{\alpha}}(\Omega)}^{2} \leq K(T_{0}).
$$
Letting $T\to \infty$ we conclude that $u^{\beta+1}\in
L^{2^*_{\alpha}}(\Omega)$.
Clearly we can obtain that $f(\cdot,u)\in L^{r}(\Omega) $
for some $r>N/\alpha$, in a finite number of steps. Thus, we conclude
applying Theorem 4.7 of \cite{brandle-colorado-depablo-sanchez}.
~\hfill$\Box$

Now we characterize the regularity of the solutions of
$(P_{\lambda})$ for the whole range of exponents.
\begin{Proposition}\label{prop:reg}
Let $u$ be a solution of $(P_{\lambda})$. Then the following hold
\begin{enumerate}
\item[(i)] If $\alpha=1$ and $q\geq 1$ then $u \in\mathcal{C}^{\infty}(\overline{\Omega}).$
\item[(ii)] If $\alpha=1$ and $q < 1$ then $u \in\mathcal{C}^{1,q}(\overline{\Omega}).$
\item[(iii)] If $\alpha< 1$ then $u \in\mathcal{C}^{\alpha}(\overline{\Omega})$.
\item[(iv)] If $\alpha>1$ then $u\in\mathcal{C}^{1,\alpha-1}(\overline{\Omega})$.
\end{enumerate}
\end{Proposition}

\pf First we observe that, by Proposition \ref{prop:bound}, we have
$u\in L^{\infty}(\Omega)$ and also $f_{\lambda}(u)\in
L^{\infty}(\Omega).$

\begin{itemize}
\item[(i)]{Applying Proposition 3.1 of
\cite{cabre-tan}, we get that $u\in \mathcal{C}^{\gamma}(\overline{\Omega})$, for some $\gamma<1$. Since
$q\geq1$ then $f_{\lambda}(u) \in
\mathcal{C}^{\gamma}(\overline{\Omega})$, so, again by Proposition 3.1 of
\cite{cabre-tan}, it follows that $u\in\mathcal{C}^{1,\gamma}(\overline{\Omega})$. Iterating the process we conclude that $u \in
\mathcal{C}^{\infty}(\overline{\Omega})$. }
\item[(ii)]{As before we have $u\in\mathcal{C}^{\gamma}(\overline{\Omega})$, for some $\gamma<1$. Therefore $f_{\lambda}(u)
\in \mathcal{C}^{q\gamma}(\overline{\Omega})$. It follows that $u\in\mathcal{C}^{1,q\gamma}(\overline{\Omega})$, which gives $f_{\lambda}(u)\in
\mathcal{C}^{q}(\overline{\Omega})$. Finally this implies $u\in\mathcal{C}^{1,q}(\overline{\Omega}).$}
\item[(iii)]{By Lemma 2.8 of \cite{capella-davila-dupaigne-sire}
we obtain that $u\in \mathcal{C}^{\gamma}(\overline{\Omega})$ for all $\gamma\in(0,\alpha)$.
This implies that $f_{\lambda}(u)\in \mathcal{C}^{r}(\overline{\Omega})$ for
every $r<\min\{q\alpha,\,\alpha\}$. Therefore, again by
\cite{capella-davila-dupaigne-sire}, this time using Lemmas 2.7 and 2.9,
we get that $u\in \mathcal{C}^{\alpha}(\overline{\Omega})$.}
\item[(iv)] Since $\alpha>1$, we can write problem $(P_{\lambda})$ as
follows
\begin{equation}\label{formulacion_equivalente}
   \left\{\begin{array}{ll}
    (-\Delta)^{1/2}u=s&\quad\mbox{in }\Omega,\\
    (-\Delta)^{(\alpha-1)/2}s=f_{\lambda}(u)&\quad\mbox{in }
    \Omega,\\
    u=s=0&\quad\mbox{on } \partial\Omega.
  \end{array}\right.
\end{equation}
Reasoning as before, we obtain the desired regularity in two steps, using
Proposition~3.1 in
\cite{cabre-tan} and Lemmas~2.7 and~2.9 in \cite{capella-davila-dupaigne-sire}.
\end{itemize}
~\hfill$\Box$

We end this section adapting  to our setting a concentration-compactness result
by P.L. Lions
\cite{Lions85b}, used in the proof of Lema \ref{strongly_convergent}. We recall
that a related concentration-compactness result for the fractional Laplacian
has been recently obtained in \cite{pp}. Nevertheless, we need the version
corresponding  to the extended problem, and it cannot be deduced from the one
in \cite{pp}.
\begin{Theorem}\label{lions}
Let $\{w_{n}\}_{n\in\mathbb{N}}$ be a weakly convergent sequence to $w$ in
$X^{\alpha}_{0}(\mathcal{C}_{\Omega})$, such that the
sequence $\{y^{1-\alpha}|\nabla w_{n}|^{2}\}_{n\in\mathbb{N}}$ is tight.  Let $u_{n}=Tr(w_{n})$
and $u=Tr(w)$. Let $\mu$, $\nu$ be two non negative measures such that
 \begin{equation}\label{hipotesis}
 y^{1-\alpha}|\nabla w_{n}|^{2} \to \mu \qquad \mbox{and} \qquad |u_{n}|^{2^{*}_{\alpha}} \to \nu,\quad \mbox{as }\, n\to \infty
 \end{equation}
in the sense of measures.
 Then there exist an at most countable set $I$ and points $\{x_{i}\}_{i\in I}\subset\Omega$ such that
\begin{enumerate}
\item  $ \displaystyle \nu = |u|^{2^{*}_{\alpha}} + \sum_{k\in I}
\nu_{k}\delta_{x_{k}}$, $\nu_{k}>0,$
\item  $ \displaystyle \mu
\geq y^{1-\alpha}|\nabla w|^{2} + \sum_{k\in I}
\mu_{k}\delta_{x_{k}}$, $\mu_{k}>0,$
\item  $
\displaystyle\mu_{k}\geq
S(\alpha,N)\nu_{k}^{\frac{2}{2^*_{\alpha}}}.$
\end{enumerate}
\end{Theorem}
\pf Let $\varphi\in
C_{0}^{\infty}(\overline{\mathcal{C}_{\Omega}})$. By the trace inequality
\eqref{eq:trace-r} with $r=2^*_\alpha$ it follows that
\begin{equation}\label{trazas}
S(\alpha,N)\left(\int_{\Omega}{|\varphi
w_{n}|^{2^{*}_{\alpha}}dx}\right)^{2/2^{*}_{\alpha}}\leq
\int_{\mathcal{C}_{\Omega}}{y^{1-\alpha}|\nabla(\varphi
w_{n})|^{2}dxdy}.
\end{equation}
Let $K^{\ast}:=K_{1}\times K_{2}\subseteq
\overline{\mathcal{C}}_{\Omega}$ be the support of $\varphi$ and
suppose first that the weak limit $w=0$. Then we get that
\begin{eqnarray*}
\int_{\mathcal{C}_{\Omega}}{y^{1-\alpha}|\nabla(\varphi w_{n})|^{2}dxdy}&=&\int_{K^{\ast}}{y^{1-\alpha}|\nabla(\varphi w_{n})|^{2}dxdy}\\
&=&\int_{K^{\ast}}{y^{1-\alpha}|w_{n}|^{2}|\nabla\varphi|^{2}dxdy}+\int_{K^{\ast}}{y^{1-\alpha}|\varphi|^{2}|\nabla w_{n}|^{2}dxdy}\\
&& + 2\int_{K^{\ast}}{y^{1-\alpha}w_{n}\varphi\langle\nabla\varphi,\nabla{w_{n}}\rangle
dxdy}.
\end{eqnarray*}
Since $K^{\ast}$ is a bounded domain, and $y^{1-\alpha}$ is an $A_2$ weight, we
have the compact inclusion
\[H^{1}(K^{\ast},y^{1-\alpha})\hookrightarrow\hookrightarrow L^{r}(K^{\ast}, y^{1-\alpha})\,,\,1\leq r < \frac{2(N+1)}{N-1}\,,\,\alpha\in (0,2).\]
Therefore, for a suitable subsequence, we get the limit
\[\int_{K^{\ast}}{y^{1-\alpha}|w_{n}|^{2}|\nabla\varphi|^{2}dxdy}\rightarrow 0,\quad \mbox{as }\, n\to \infty.
\]
By the weak convergence, given by hypothesis, we obtain
\[\int_{K^{\ast}}{y^{1-\alpha}w_{n}\varphi\langle\nabla\varphi,\nabla{w_{n}}\rangle dxdy}\rightarrow 0,\quad \mbox{as }\, n\to \infty. \]
Hence, by (\ref{hipotesis}) we conclude that
\[\int_{\mathcal{C}_{\Omega}}{y^{1-\alpha}|\nabla(\varphi w_{n})|^{2}dxdy}\rightarrow\int_{\mathcal{C}_{\Omega}}{|\varphi(x,y)|^{2}d\mu},\quad \mbox{as }\, n\to \infty. \]
Then, from (\ref{trazas}) we get
\begin{equation}\label{RH}
S(\alpha,
N)\left(\int_{\Omega}{|\varphi|^{2^{*}_{\alpha}}d\nu}\right)^{2/2^{*}_{\alpha}}\leq
\int_{\mathcal{C}_{\Omega}}{|\varphi|^{2}d\mu},\quad \forall\;\varphi\in
C_{0}^{\infty}(\overline{\mathcal{C}}_{\Omega}).
\end{equation}
If now $w\neq 0$, we apply the above result to the function $v_{n}=w_{n}-w$.
Indeed if
\[y^{1-\alpha}|\nabla v_{n}|^{2} \to d\tilde{\mu} \qquad \mbox{and} \qquad |v_{n}(\cdot,0)|^{2^{*}_{\alpha}} \to d\tilde{\nu},\quad \mbox{as }\, n\to \infty,\]
it follows that
\[S(\alpha, N)\left(\int_{\Omega}{|\varphi|^{2^{*}_{\alpha}}
d\tilde{\nu}}\right)^{2/2^{*}_{\alpha}}\leq\int_{\mathcal{C}_{\Omega}}{|\varphi|^{2}d\tilde{\mu}}\,,
\quad \forall\;\varphi\in C_{0}^{\infty}(\overline{\mathcal{C}}_{\Omega}),\] therefore,
(\cite{Lions85b}), for some sequence of points $\{x_{k}\}_{k\in I}\subset\Omega$, we have
\[d\tilde{\nu}=\sum_{k\in I}{\tilde{\nu}_{k}\delta_{x_{k}}}\,,\qquad
d\tilde{\mu}\geq\sum_{k\in I}{\tilde{\mu}_{k}\delta_{x_{k}}}\,,\] with
$\tilde{\mu}_{k}\geq S(\alpha,N)\tilde{\nu}_{k}^{2^*_\alpha/2}$. Hence, by
Lemma
\ref{brezis_lieb}, we obtain
$$
d\nu=|u|^{2^{*}_{\alpha}}+\sum_{k\in I}{\tilde{\nu}_{k}\delta_{x_{k}}}.$$
Let now $\varphi$ be a test function. We have
$$
\begin{array}{rcl}
\displaystyle\int_{\mathcal{C}_{\Omega}}{y^{1-\alpha}\varphi|\nabla w_{n}|^{2}dxdy}&=&\displaystyle\int_{\mathcal{C}_{\Omega}}{y^{1-\alpha}\varphi|\nabla w|^{2}dxdy}+\int_{\mathcal{C}_{\Omega}}{y^{1-\alpha}\varphi|\nabla (w_{n}-w)|^{2}dxdy}\\ & & \\
& &\displaystyle +2\int_{\mathcal{C}_{\Omega}}{y^{1-\alpha}\varphi\langle\nabla
(w_{n}-w),\nabla w\rangle dxdy}.
\end{array}
$$
Taking limits as $n\to\infty$ we get that
$$
\begin{array}{rcl}
\displaystyle\int_{\mathcal{C}_{\Omega}}{\varphi d\mu}&=&\displaystyle\int_{\mathcal{C}_{\Omega}}{y^{1-\alpha}\varphi|\nabla w|^{2}dxdy}+\int_{\mathcal{C}_{\Omega}}{\varphi d\tilde{\mu}}\\ & &\\
&\geq&\displaystyle\int_{\mathcal{C}_{\Omega}}{y^{1-\alpha}\varphi|\nabla
w|^{2}dxdy}+\int_{\mathcal{C}_{\Omega}}{y^{1-\alpha}\varphi\sum_{k\in
I}{\tilde{\mu}_{k}\delta_{x_{k}}}dxdy},
\end{array}
$$
with the same condition $\tilde{\mu}_{k}\geq
S(\alpha,N)\tilde{\nu}_{k}^{2^*_\alpha/2}.$ So we obtain the desired
conclusion.~\hfill$\Box$

\

\bf{Acknowledgements}:
\rm{The authors want to thank professor Fernando Soria for his valuable
help about the Poisson kernel estimates. We also  thank professor Jes\'{u}s
Garc{\'\i}a-Azorero for some helpful discussions.}


\end{document}